\documentclass[12pt]{article}
\usepackage{amssymb}
\usepackage{amsfonts}
\usepackage{mathrsfs}
\usepackage{graphicx}
\usepackage{amsmath}
\usepackage{epsfig}
\usepackage{verbatim}
\usepackage{url}

    \newcommand{\lab}[1]{\label{#1}}                

\def\smallpage{
\addtolength\textwidth{3cm} \addtolength\oddsidemargin{-1.5cm}
\addtolength\textheight{3cm} \addtolength\topmargin{-1.5cm}}
 \smallpage

\newcommand{\remove}[1]{}
\newcommand{\tremove}[1]{}
\newcommand{\bel}[1]{\be\lab{#1}}
\newcommand\eqn[1]{(\ref{#1})}

\newcommand{\be}{\begin{equation}}
\newcommand{\ee}{\end{equation}}
\newcommand{\bea}{\begin{eqnarray}}
\newcommand{\eea}{\end{eqnarray}}
\newcommand{\bean}{\begin{eqnarray*}}
\newcommand{\eean}{\end{eqnarray*}}

\newtheorem{thm}{Theorem}[section]
\newtheorem{cor}[thm]{Corollary}

\newtheorem{lemma}[thm]{Lemma}
\newtheorem{definition}[thm]{Definition}
\newtheorem{claim}[thm]{Claim}
\newtheorem{prop}[thm]{Proposition}
\def\proof{\noindent{\bf Proof.\ }  }
\def\qed{~~\vrule height8pt width4pt depth0pt}

\def\ex{{\bf E}}
\def\pr{{\bf P}}

\def\hD{{\hat \D}}

\def\coreH{\widehat H}
\def\core2H{G}
\def\orH{H}
\def\barS{\overline S}
\def\barA{\overline A}
\def\Ll{Q}
\def\D{{\mathcal D}}
\def\hD{D_0}

\def\eps{\epsilon}
\def\la{\lambda}

\def\ss {\smallskip}

\def\no{\noindent}

\catcode`@=11 \@addtoreset{equation}{section}

\catcode`@=12

\title{Orientability thresholds for random hypergraphs}
\author{Pu Gao\footnote{The research reported in this paper
forms part of this author's Ph.D. thesis ``Generation and properties
of random graphs and analysis of randomized
  algorithms'' submitted to the
University of Waterloo in 2009.}\ \ \footnote{Current affiliation:
Max-Planck-Institut f\"{u}r Informatik, 66123 Saarbr\"{u}cken, Saarland, Germany}\\
 Department of Combinatorics and Optimization\\
University of Waterloo, Canada
 \\ janegao@mpi-inf.mpg.de,
\and Nicholas Wormald\thanks{Research
supported by the Canadian Research Chairs Program and NSERC}\\ Department of Combinatorics and Optimization\\
University of Waterloo, Canada\\ nwormald@math.uwaterloo.ca}
\date{}

\begin{document}

\maketitle
\begin{abstract}

Let $h>w>0$ be two fixed integers. Let $\orH$ be a random hypergraph
whose hyperedges are all of cardinality $h$.  To {\em $w$-orient} a
hyperedge, we assign exactly $w$ of its vertices positive signs with
respect to the hyperedge, and the rest negative.  A
$(w,k)$-orientation of $\orH$ consists of a $w$-orientation of all
hyperedges of $\orH$, such that each vertex receives at most $k$
positive signs from its incident hyperedges.  When $k$ is large
enough, we determine the threshold of the existence of a
$(w,k)$-orientation of a random hypergraph.
   The $(w,k)$-orientation
of hypergraphs is strongly related to a general version of the
off-line load balancing problem. The graph case, when $h=2$ and
$w=1$, was solved recently by Cain, Sanders and Wormald and
independently by Fernholz and Ramachandran, which settled a
conjecture of Karp and Saks.
\end{abstract}

 \section{Introduction}
\lab{introduction}

In this paper we consider a generalisation to random hypergraphs of
a commonly studied orientation problem on graphs.
 An {\em $h$-hypergraph} is a hypergraph whose hyperedges are all of
size $h$. Let $h>w$ be two given positive integers. We consider
$\mathcal{G}_{n, m,h}$, the probability space of the set of all
$h$-hypergraphs on $n$ vertices and $m$ hyperedges with the uniform
distribution. A hyperedge is said to be {\em $w$-oriented} if
exactly $w$ distinct vertices in it are marked with positive signs
with respect to the hyperedge. The {\em indegree} of a vertex is the
number of positive signs it receives. Let $k$ be a positive integer.
A $(w,k)$-orientation of an $h$-hypergraph is a $w$-orientation all
hyperedges such that each vertex has indegree at most $k$. If such a
$(w,k)$-orientation exists, we say the hypergraph is {\em
$(w,k)$-orientable}; for $w=1$ we simply say {\em $k$-orientable.}
Of course, being able to determine the
 $(w,k)$-orientability of an $h$-hypergraph
$\orH$ for all $k$ solves the optimisation problem of minimising the
maximum indegree of a $w$-orientation of $\orH$. If a graph (i.e.\
the case $h=2$) is $(1,k)$-oriented, we may orient each edge of the
graph in the normal fashion towards its vertex of positive sign, and
we say the graph is {\em $k$-oriented}.

Note that a sufficiently sparse hypergraph is easily
$(w,k)$-orientable. On the other hand, a trivial requirement for
$(w,k)$-orientability is $m\le kn/w$, since any $w$-oriented
$h$-hypergraph with $m$ edges has average indegree $mw/n$. In this
paper, we show the existence and determine the value of the sharp
threshold (defined more precisely later) at which the random
$h$-hypergraph $\mathcal{G}_{n,m,h}$ fails to be $(w,k)$-orientable,
provided $k$ is a sufficiently large constant. \remove{ To be
precise, this threshold is a number $c_{h,w,k}$ such that
$\mathcal{G}_{n,m,h}$ is asymptotically almost surely (a.a.s.)
$(w,k)$-orientable for $m<(c_{h,w,k}-\eps)n$, and a.a.s.\ {\em not}
$(w,k)$-orientable for $m>(c_{h,w,k}+\eps)n$, for any fixed
$\eps>0$. }  We show that the threshold is the same as the threshold
at which a certain type of subhypergraph achieves a critical
density. In the above, as elsewhere in this paper, the phrase ``for
k sufficiently large'' means  for $k$ larger than some constant
depending only on  $w$ and $h$.

\remove{
\subsection{Applications to load balancing, and\nn{} previous results}
}
The hypergraph orientation problem is motivated by classical load
balancing problems which have appeared in various guises in computer
networking.  A seminal result of Azar, Broder, Karlin and
Upfal~\cite{ABKU} is as follows. Throw $n$ balls
  sequentially into $n$ bins, with each ball put into the least-full of $h\ge 2$
   randomly chosen boxes. Then, with high probability,  by the time all balls are allocated,
    no bin contains many more than $(\ln \ln n) /\ln h$ balls.  If, instead,
    each ball is placed in a random bin, a much larger maximum value is likely
    to occur, approximately $\ln n /(\ln \ln n)$.
This surprisingly simple method of reducing the maximum is widely
used for load balancing. It has become known as the multiple-choice
paradigm, the most common version being two-choice, when $h=2$.

One application of load balancing occurs when work is spread among a
group of computers, hard drives, CPUs, or other resources. In the
on-line version, the jobs arrive sequentially and are assigned to
separate machines. To save time, the load balancer decides which
machine a job goes to after checking the current load of only a few
(say $h$) machines. The goal is to minimise the maximum load   of a
machine.  Mitzenmacher,  Richa and Sitaraman~\cite{MRS}  survey the
history, applications and techniques relating to this. In
particular,    Berenbrink,   Czumaj,   Steger, and
V\H{o}cking~\cite{ABKU,BCSV}   show an achieveable  maximum load is
$m/n+O(\log \log n)$ for $m$ jobs and $n$ machines when $h\ge 2$.

 Another application of load balancing, more relevant to the topic of this paper, is mentioned by Cain, Sanders and
 the second author~\cite{CSW}. This is the disk scheduling problem, in the context where any
$w$ out of $h$ pieces of data are needed to reconstruct a logical
data block. Individual pieces can be initially stored on different
disks. Such an arrangement has advantageous   fault tolerance
features to guard against disk failures. It is also good for load
balancing: when a request for a data block arrives, the scheduler
can choose  any $w$ disks among the $h$ relevant ones. See Sanders,
Egner  and Korst~\cite{SEK} for further references.

These load balancing problems correspond to the $(w,k)$-orientation
problem for $h$-uniform hypergraphs, with $w=1$ in the case of the
job scheduling problem. The machines (bins) are vertices and a job
(ball)  is an edge consisting of precisely the set of machines to
which it can be allocated. A job is allocated to a machine by
assigning a positive sign to that vertex. The maximum load is then
equal to the maximum indegree of a vertex in the $(w,k)$-oriented
hypergraph.

The work in this paper is motivated by the off-line version of this
problem, in which the edges are all exposed at the start. This has
obvious applications, for instance, in the disk scheduling problem,
the scheduler may be able to quickly process a large number of requests together off-line, to balance the load better.
 This can be useful if there is a backlog of requests; of course, if backlogs do not occur, the online problem is more relevant,
 but this would indicate ample processing capacity, in which case there may be less need for load balancing in the first place.
 Trivially, the on-line and off-line versions are the same if $h=1$, i.e.\ there is no choice.  For $h=2$, the off-line version experiences an even better improvement than the on-line one.  If $m<cn$ items are allocated to $n$ bins, for $c$ constant,  the expected maximum load is bounded above by some constant $c'$ depending on $c$.

 To our knowledge, previous theoretical results concern only the case
 $w=1$ (this also applies to on-line). For $w=1$
it is well known that an optimal off-line solution, i.e.\ achieving
minimum possible maximum load, can be found in polynomial time
($O(m^2)$) by solving a maximum flow problem.  As explained
in~\cite{CSW}, it is desirable to achieve fast algorithms that are
close to optimal with respect to the maximum load. There are linear
time algorithms that achieve maximum load
$O(m/n)$~\cite{DM,KLM,MSS}.

A central role  in solutions of the off-line  orientation problem
 with $(w,h)=(1,2)$ is played by the  {\em $k$-core}
of a graph, being the largest subgraph with minimum degree at least
$k$. The sharp threshold   for the $k$-orientability of the random
graph $\mathcal{G}(n,m)= \mathcal{G}_{n,m,2}$ was found
in~\cite{CSW}, and simultaneously by Fernholz and
Ramachandran~\cite{FR}.    These were proofs of a conjecture of Karp
and Saks  that this threshold coincides with the threshold at which
the $(k+1)$-core has average degree at most $2k$. (It is obvious
that a graph cannot be $k$-oriented if it has a subgraph of average
degree greater than $2k$.)  In each case, the proof
  analysed a linear time algorithm that finds
a $k$-orientation a.a.s.\ when the mean degree of the $(k+1)$-core
is slightly less than $2k$. In this sense, the algorithms are
asymptotically optimal since the threshold for the algorithms
succeeding coincides with the threshold for existence of a
$k$-orientation. The proof in~\cite{FR} was significantly simpler
than the other, which was made possible because a different
algorithm was employed. It used a trick of ``splitting vertices'' to
postpone decisions and thereby reduced the number of variables to be
considered.

During the preparation of this paper, three preprints appeared by
Frieze and  Melsted~\cite{FM}, Fountoulakis and
Panagiotou~\cite{FP}, and by Dietzfelbinger, Goerdt, Mitzenmacher,
Montanari, Pagh and Rink~\cite{DGMMPR} which independently study the
threshold of $(1,1)$-orientability of $\mathcal{G}_{n,m,h}$, i.e.\
the case $w=k=1$. This has applications to cuckoo hashing.  However,
there seems to be no easy way to extend the proofs
in~\cite{DGMMPR,FM,FP} to solve for the case $k>1$, even when $w=1$.

In this paper, we solve the generalisation of the conjecture of Karp
and Saks mentioned above, for fixed  $h>w>0$, provided   $k$ is
sufficiently large. That is, we find   the threshold of
$(w,k)$-orientation of random $h$-hypergraphs in $\mathcal{G}_{n,
m,h}$. The determination of this threshold   helps to predict loads
in the off-line $w$-out-of-$h$ disk scheduling problem, where the
randomness of the hypergraph is justified by the random intial
allocation of file segments to disks.
 We believe
furthermore that the characterisation of the threshold in terms of
density of a type of core, and possibly our method of proof, will
potentially help lead to fast algorithms for finding asymptotically
optimal orientations.

Our approach has a significant difference from that used in the
graph case when $(w,h)=(1,2)$. The algorithm used in~\cite{FR} does
not seem to apply in the hypergraph case, at least, splitting
vertices cannot be done without creating hyperedges of larger and
larger size. The algorithm used by~\cite{CSW}, on the other hand,
generalises in an obvious way, but it is already very complicated to
analyse  in the graph case, and the extension of the analysis to the
hypergraph case seems formidable. However, in common with those two
approaches, we first find what we call the $(w,k+1)$-core in the
hypergraph, which is an analogue of the $(k+1)$-core in graphs. We
 determine the size and
density of this core when the random hypergraph's density is
significantly larger than what is required for the core to form.
This result may be of independent interest, and uses the differential equation method in a setting
which contains a twist not encountered when it is applied to the graph case: some of the functions involved have singularities at the starting point. (See Section~\ref{core} for details.)

 Although we
gain information on the threshold of appearance of this core, we do
not, and do not need to, determine it precisely. From here,   we use
the natural representation of the orientation problem in terms of
flows.
 It is quite easy to generalise the network flow formulation from the
 case $h=2$, $w=1$ to the arbitrary case, giving a problem that can be
 solved in time $O(m^2)$ for  $m=\Theta(n)$. Unlike the approaches for
 the graph case, we do not study an algorithm that solves the load balancing
 problem. Instead, we use the minimum cut characterisation of the maximum flow
 to  show that a.a.s.\ the hypergraph   can be $(w,k)$-oriented if and only
 if the density of its $(w,k+1)$-core is below a certain threshold. When the
 density of the  $(w,k+1)$-core is above this threshold, it is trivially too dense to be $(w,k)$-oriented.
  Even the case $w=1$ of
our result gives a significant generalisation of the known results.
We prove that the threshold of the orientability coincides with the
threshold at which certain type of density (in the case $w=1$, this
refers to the average degree divided by $h$) of the $(w,k+1)$-core
is at most $k$, and also the threshold at which certain type of
induced subgraph (in the case $w=1$, this refers to the standard
induced subgraph) does not appear. For  the graph case, our method
 provides a new proof
(for sufficiently large $k$) of the Karp-Saks conjecture that we
believe is simpler than the proofs of~\cite{CSW} and~\cite{FR}.

 We give precise statements of our results, including definition of the $(w,k+1)$-core,
 in Section~\ref{mainResult}.  In Section~\ref{core}  we study
 the properties of the
 $(w,k+1)$-core.
 In Section~\ref{orientability}, we formulate the appropriate network flow problem,
 determine a canonical minimum cut for a network corresponding to a non-$(w,k)$-orientable hypergraph,
 and give conditions under which such a minimum cut can exist.
Finally, in Section~\ref{probabilistic}, we show that for $k$ is
sufficiently large, such a cut  a.a.s.\ does not exist when the
density of the core is below a certain threshold.

An extended abstract for this paper, omitting most proofs, will
appear in STOC 2010~\cite{GW4}.

\section{Main results}
\lab{mainResult}

Let $h>w>0$ and $k\ge 2$ be fixed. For any $h$-hypergraph $\orH$, we
examine whether a $(w,k)$-orientation exists. We call a vertex  {\em
light} if the degree of the vertex is at most $k$. For any light
vertex $v$, we can give $v$ the positive sign respect to any
hyperedge $x$ that is incident to $v$ (we call this {\em partially
orienting} $x$ towards to $v$), without violating the condition that
each vertex has indegree at most $k$. Remove $v$ from $\orH$, and
for each hyperedge $x$ incident to $v$, simply update $x$ by
removing $v$. Then the size of $x$ decreases by $1$, and it
has one less vertex that needs to be given the positive sign. 
If the size of a hyperedge falls to $h-w$,
 we can simply remove that hyperedge from the hypergraph. Repeating this until no light vertex
exists, we call the remaining hypergraph $\coreH$ the {\em
$(w,k+1)$-core} of $\orH$. Every vertex in $\coreH$ has degree at
least $k+1$, and every hyperedge in $\coreH$ of size $h-j$ requires
a $(w-j)$-orientation in order to obtain a $w$-orientation of the
original hyperedge in $\orH$.

In order to simplify the notation, we use $\bar n$, $\bar m$ and
$\bar \mu$ to denote the numbers of vertices and of hyperedges, and
the average degree, of $H\in\mathcal{G}_{\bar n,\bar m,h}$,
reserving $n$, $m_{h-j}$ and $\mu$ to denote the numbers of vertices
and of hyperedges of size $h-j$, and the average degree, of
$\coreH$.

Instead of considering the probability space $\mathcal{G}_{\bar
n,\bar m,h}$, we may consider $\mathcal{M}_{\bar n,\bar m,h}$, the
probability space of random multihypergraphs with $\bar n$ vertices
and $\bar m$ hyperedges, such that each hyperedge $x$ is of size
$h$, and each vertex in $x$ is chosen independently, uniformly at
random from $[\bar n]$. Actually $\mathcal{M}_{\bar n,\bar m,h}$ may
be a more accurate model for the off-line load balancing problem in
some applications, and as we shall see, results for the non-multiple
edge case can be deduced from it. For a nonnegative integer vector
${\bf m}=(m_2,\ldots,m_h)$, we also define the probability space
$\mathcal{M}_{n,{\bf m}}$, being the obvious generalisation of
$\mathcal{M}_{n,m,h}$ to  non-uniform multihypergraphs  in which
$m_i$ is the number of hyperedges of size $i$.

 All our asymptotic notation refers to $n\to\infty$.  For clarity, we consider $H\in\mathcal{M}_{\bar n,\bar m,h}$. We use
$n$, $m_{h-j}$ and $\mu$ for the number of vertices, the number of
hyperedges of size $h-j$ and the average degree of $\coreH$. We
parametrise the number $\bar m$ of edges in the hypergraphs under
study by letting
  $\bar \mu = \bar\mu(n)$ denote $h\bar m/\bar n$, the average degree of
$H\in\mathcal{M}_{\bar n,\bar m,h}$ (or  of $H\in  \mathcal{G}_{\bar
n,\bar m,h}$).

 Our first observation
concerns the distribution of $\coreH$ and its vertex degrees.
 Let ${\rm Multi}(n,m,k+1)$ denote the
  multinomial  distribution of  $n$ integers summing to $m$, restricted
to each of the  integers being at least $k+1$.
  We call this the {\em
truncated multinomial} distribution.
\begin{prop}\lab{uniformCore}
Let $h>w\ge 1$ be two fixed integers. Let $\orH\in\mathcal{M}_{\bar
n,\bar m,h}$ and let $\coreH$ be its $(w,k+1)$-core. Conditional on
its number $n$ of vertices and numbers $m_{h-j}$ of hyperedges of
size $h-j$ for $j=0,\ldots,w-1$, the random hypergraph $\coreH$ is
distributed uniformly at random.\remove{randomly on multihypergraphs
with all vertices of degree at least $k+1$ and having the \nn{}given
parameters $n$ and $m_{h-j}$ for each $j$.} Furthermore, the
distribution of the degree sequence of $\coreH$ is  the truncated
multinomial distribution ${\rm Multi}(n,m,k+1)$ where
$m=\sum_{j=0}^{w-1} (h-j)m_{h-j}$.
\end{prop}

 The following theorem shows that the size and the number of hyperedges of
 $\coreH$ are highly concentrated around the solution of a system of
differential equations. The theorem covers the cases for any
arbitrary $h>w\ge 2$ and holds for all sufficiently large $k$. The
special case $w=1$ has been studied by various authors and the
concentration results can be found in~\cite[Theorem 3]{CSW} which
hold for all $k\ge 0$. Since $w$ and $k$ are fixed, we often omit
them from the notation.
\begin{thm}\lab{t:CoreDensity}
Let $h>w\ge 2$ be two fixed integers. Assume that  for some constant
$c>1$ we have $ck\le \bar\mu$ where $\bar \mu= h\bar m/\bar n$. Let
$\orH\in\mathcal{M}_{\bar n,\bar m,h}$ and let $\coreH$ be its
$(w,k+1)$-core. Let $n$ be the number of vertices  and $m_{h-j}$ the
number of hyperedges of size $h-j$ of $\coreH$.  Then, provided $k$
is sufficiently large, there are constants $\alpha>0$ and
$\beta_{h-j}>0$, defined in~\eqn{alphabeta} below, depending  only
on $\bar \mu$, $k$, $w$ and $h$, for which a.a.s.\ $n\sim \alpha\bar
n$ and $m_{h-j}\sim \beta_{h-j}\bar n$
 for $0\le i\le w-1$.  The same conclusion (with the same constants) holds for
$\orH\in\mathcal{G}_{\bar n,\bar m,h}$.
\end{thm}
\noindent
{\bf Note.\ } The full definition of $\alpha$ and $\beta_{h-j}$ in the theorem  is rather complicated, involving the solution of a differential equation system  given below in~(\ref{eq:DL}--\ref{eq:DefSystem}).

Let $\mathcal{P}$ be a hypergraph property and let
$\mathcal{M}_{n,m,h}\in\mathcal{P}$ denote the event that a random
hypergraph from $\mathcal{M}_{n,m,h}$ has the property
$\mathcal{P}$. Following~\cite[Section 10.1, Definition 4]{AS}, we
say that $\mathcal{P}$ has a {\em sharp threshold function} $f(n)$
if for any constant $\eps>0$,
$\pr(\mathcal{M}_{n,m,h}\in\mathcal{P})\rightarrow 1$ when
$m\le(1-\eps) f(n)$, and
$\pr(\mathcal{M}_{n,m,h}\in\mathcal{P})\rightarrow 0$ when
$m\ge(1+\eps) f(n)$.

  Let
$\kappa(\coreH)$ denote $\sum_{j=0}^{w-1}(w-j)m_{h-j}/n$, which we
call the {\em $w$-density} of $\coreH$. We similarly define the
$w$-density of any hypergraph all of whose hyperedges have sizes
between $h-w+1$ and $h$. It helps to notice, by the definition of
$w$-density, that
$$n\kappa(\coreH)=d(\coreH)-(h-w)m,
$$  where $d(\coreH)$ denotes the
degree sum of $\coreH$ and $m=\sum_{j=0}^{w-1}m_{h-j}$. We say that
a hypergraph $\orH$ has {\em property} $\mathcal{T}$ if
$\kappa(\coreH)\le k$, where $\coreH$ is the $(w,k+1)$-core of $H$.
 The
following theorem, proved using Theorem~\ref{t:CoreDensity},
immediately gives the corollary that there is a sharp threshold
function for property $\mathcal{T}$.
\begin{thm}\lab{t2:CoreDensity} Let $\orH\in\mathcal{M}_{\bar n,\bar
m,h}$. Let $\bar\mu$ be the average degree of $\orH$ and let
$\coreH$ be the $(w,k+1)$-core of $\orH$. Then for all sufficiently
large $k$, there exists a strictly increasing function $c(\bar\mu)$
of $\bar\mu$, such that for any fixed $c_2>c_1>1$ and for any
$c_1k<\bar\mu<c_2k$, a.a.s.\ $\kappa(\coreH)\sim c(\bar\mu)$.
\end{thm}
\begin{cor}\lab{c:CoreDensity}
There exists a sharp threshold function $f(\bar n)$ for the
hypergraph property $\mathcal{T}$ in $\mathcal{M}_{\bar n,\bar m,h}$
and $\mathcal{G}_{\bar n,\bar m,h}$ provided $k$ is sufficiently
large.
\end{cor}
 The function $c(\bar\mu)$ in the theorem, and the
threshold function in the corollary, are determined by the solution
of the differential equation system referred to in Theorem 2.1.

We have defined a $(w,k)$-orientation of a uniform hypergraph in
Section~\ref{introduction}. We can
 similarly define a
$(w,k)$-orientation of a non-uniform hypergraph $\core2H$ with sizes
of hyperedges between $h-w+1$ and $h$ to be a simultaneous
$(w-j)$-orientation of each hyperedge of size $h-j$ such that every
vertex has indegree at most $k$.
 By counting the positive signs in orientations, we see that if
property  $\mathcal{T}$ fails, there is no $(w,k)$-orientation of
$\coreH$, and hence there is no $(w,k)$-orientation of $\orH$.

  For a nonnegative integer vector ${\bf m}=(m_{h-w+1},\ldots,m_h)$,
 let $\mathcal{M}(n, {\bf m}, k+1)$ denote $\mathcal{M}_{n,{\bf m}}$
restricted to multihypergraphs with minimum degree at least $k+1$.
By Proposition~\ref{uniformCore}, $\mathcal{M}(n,{\bf m},k+1)$ has
the distribution of the $(w,k+1)$-core of $H\in \mathcal{M}_{\bar
n,\bar m,h}$ conditioned on the number of vertices being $n$ and the
number of hyperedges of each size being given by ${\bf m}$. To
emphasise the difference, we will use $\core2H$ to denote a
not-necessarily-uniform hypergraph in cases where we might use
$\orH$ for a uniform hypergraph.

Given a vertex set $S$, we say a hyperedge $x$ is {\em partially
contained in $S$} if $|x\cap S|\ge 2$.
\begin{definition}\lab{d:A}
Let $0<\gamma<1$. We say that a multihypergraph $\core2H$ has
property $\mathcal{A}(\gamma)$ if for all $S\subset V(\core2H)$ with
$|S|<\gamma |V(\core2H)|$ the number of hyperedges partially
 contained in $S$ is strictly less than $k|S|/2w$.
\end{definition}

In the following theorem, ${\bf m}={\bf m}(n)$  denotes an integer
vector for each $n$.
\begin{thm}\lab{t:ForCore}
Let $\gamma$ be any constant between $0$ and $1$. Then there exists
a constant $N>0$ depending only on $\gamma$, such that for all $k>N$
and any $\eps>0$, if ${\bf m}(n)$ satisfies
$\sum_{j=0}^{w-1}(w-j)m_{h-j}(n)\le kn-\eps n$ for all $n$, then
  $\core2H\in\mathcal{M}(n,{\bf m}(n),k+1)$ a.a.s.\ either
has a $(w,k)$-orientation or  does not have property
$\mathcal{A}(\gamma)$.
\end{thm}

Let $f(\bar n)$ be the threshold of property $\mathcal{T}$ given in
Corollary~\ref{c:CoreDensity}. We show in the forthcoming
Corollary~\ref{c:SmallSet2} that for certain values of $\gamma$,
 a.a.s.\ $\coreH$ has
property $\mathcal{A}(\gamma)$ if the average degree of $\orH$ is at
most $hk/w$. We will combine this with Corollary~\ref{c:CoreDensity}
and Theorem~\ref{t:ForCore} and a relation we will show between
$\mathcal{M}_{\bar n,\bar m,h}$ and $\mathcal{G}_{\bar n,\bar m,h}$
(Lemma~\ref{uniformSimple}), to obtain the following.
\begin{cor}\lab{c2:lbalancing}
Let $h>w>0$ be two given integers and $k$ be a sufficiently large
constant. Let $f(\bar n)$ be the threshold function of property
$\mathcal{T}$ whose existence is asserted in
Corollary~\ref{c:CoreDensity}. Then $f(\bar n)$ is a sharp threshold
for the $(w,k)$-orientability of $\mathcal{M}_{\bar n,\bar m,h}$ and
$\mathcal{G}_{\bar n,\bar m,h}$.
\end{cor}

For any vertex set $S\subset V(\orH)$, define the subgraph {\em
$w$-induced} by $S$ to be the subgraph of $\core2H$ on vertex set
$S$ with the set of  hyperedges $\{x'=x\cap S:\
 x\in \orH,\ s.t.\ |x'|\ge h-w+1\}$. Call this hypergraph
$\orH_S$. \tremove{This definition no longer needs to be here: Let
$d(\orH_S)$ denote the degree sum of vertices in the hypergraph
$\orH_S$ and let $e(\orH_S)$ denote the number of hyperedges in
$\orH_S$. However, we need to have them later at an appropriate
place, OR move them back to the definition of $\kappa$. }
  It helps to notice that $\coreH$ is the largest
$w$-induced subgraph of $\orH$ with minimum degree at least $k+1$.
{From} the above results and a relation we will show between
$\mathcal{M}_{\bar n,\bar m,h}$ and $\mathcal{G}_{\bar n,\bar m,h}$,
we will obtain the following. 
\begin{cor}\lab{c:lbalancing}
The following three graph properties of $\orH\in\mathcal{M}_{\bar
n,\bar m,h}$ (or $\mathcal{G}_{\bar n,\bar m,h}$) have the same
sharp threshold.
\begin{description}
\item{(i)} $\orH$ is $(w,k)$-orientable.

\item{(ii)} $\orH$ has property $\mathcal{T}$.

\item{(iii)} There exists  no $w$-induced subgraph $H'\subset\orH$  with $\kappa(H')\ge k$.

\end{description}

\end{cor}


\section{Analysing the size and density of the $(w,k+1)$-core}
\lab{core}

A model of generating random graphs via multigraphs, used by
Bollob\'{a}s and Frieze~\cite{BF} and Chvat\'{a}l~\cite{C2}, is
described as follows. Let $\mathcal{P}_{\bar n,\bar m}$ be the
probability space of functions $g:[\bar m]\times [2]\rightarrow
[\bar n]$ with the uniform distribution.   Equivalently,
$\mathcal{P}_{\bar n,\bar m}$ can be   described as the uniform
probability space of allocations of $2\bar m$ balls into $\bar n$
bins. A probability space of random multigraphs can be obtained by
taking $\{g(i,1), g(i,2)\}$ as an edge for each $i$. This model can
easily be extended to generate non-uniform random multihypergraphs
by letting ${\bf m}=(m_2,\ldots,m_h)$ and taking $\mathcal{P}_{\bar
n,{\bf m}}=\{g:\cup_{i=2}^h[m_i]\times [i]\rightarrow [\bar n]\}$.
Let $\mathcal{M}_{\bar n,{\bf m}}$ be the probability space of
random multihypergraphs obtained by taking each
$\{g(j,1),\ldots,g(j,i)\}$ as a hyperedge, where $j\in [m_i]$ and
$2\le i\le h$. (Loops and multiple edges are possible.) Note that
$\mathcal{M}_{\bar n,{\bf m}}$,  where ${\bf m} = (m_2)=(\bar m)$,
is a random multigraph; it was shown in~\cite{C2} that if this is
conditioned on being simple (i.e.\ no loops and no multiple edges),
it is equal to $\mathcal{G}_{\bar n,\bar m,2}$, and that the
probability of a multigraph in $\mathcal{M}_{\bar n,(\bar m)}$ being
simple is $\Omega(1)$ if $\bar m=O(\bar n)$. This result is easily
extended to the following result, using the same method
 of proof.
\begin{lemma}\lab{uniformSimple} Assume $h\ge 2$ is a fixed
integer and ${\bf m}=(m_2,\ldots,m_h)$ is a non-negative integer
vector. Assume further that $\sum_{i=2}^h m_i=O(\bar n)$. Then the
probability that a hypergraph in $\mathcal{M}_{\bar n,{\bf m}}$ is
simple is $\Omega(1)$.
\end{lemma}


Cain and Wormald~\cite{CW} recently introduced a related model to
analyse the $k$-core of a random (multi)graph or (multi)hypergraph,
including its size and degree distribution. This model is called the
{\em pairing-allocation} model.  The {\em partition-allocation}
model, as defined below, is a generalisation of the
pairing-allocation model, and analyses cores of multihypergraphs
with given numbers of hyperedges of various sizes. We will use  this
model to prove Theorem~\ref{t:ForCore} and to analyse a randomized
algorithm called the RanCore algorithm, defined later in this
section, which outputs the $(w,k+1)$-core of an input
$h$-hypergraph.

 Given
$h\ge 2$, $n$, ${\bf m}=(m_2,\ldots,m_h)$, ${\bf
L}=(l_2,\ldots,l_h)$ and a nonnegative integer $k$ such that
$D-\ell\ge kn$, where $D=\sum_{i=2}^{h} im_i$ and $\ell=\sum_{i=2}^h
l_i$, let $V$ be a set of $n$ bins, and ${\bf M}$ a collection of
pairwise disjoint sets $\{M_1,\ldots, M_h\}$, where $M_i$ is a set
of $im_i$ balls partitioned into parts, each of size $i$, for all
$2\le i\le h$. Let $\Ll$  be an additional bin to $V$. It may assist
the reader to know that $\Ll$ `represents' all the hyperedge
incidences at vertices of degree less than $k$, and $l_i$ is the
number of these incidences in edges of size $i$. The {\em
partition-allocation} model $\mathcal{P}(V,{\bf M},{\bf L},k)$ is
the probability space
of ways of allocating balls to bins in the following way. Let
$\mathcal{C}=\{c_2,\ldots,c_h\}$ be a set of colours. Colour balls
in $M_i$ with $c_i$. (The function of the colours is only to denote
the size of the part a ball lies in.)
  Then allocate the $D$ balls
uniformly at random  (u.a.r.)\  into the bins in $V\cup \{\Ll\}$ , such
that the following constraints are satisfied:
\begin{description}
\item{(i)} $\Ll$ contains exactly $\ell$ balls;
\item{(ii)} each bin in $V$ contains at least $k$ balls;
\item{(iii)}  for any $2\le i\le h$, the number of balls with colour $c_i$ that are
contained in $\Ll$ is $l_i$.
\end{description}

We call $\Ll$ the  {\em light} bin and all bins in $V$ {\em heavy}.
 To assist  with the analysis in some
situations, we consider the following  algorithm  which clearly
generates a probability space equivalent to $\mathcal{P}(V,{\bf
M},{\bf L},k)$.  We call this alternative the {\em
allocation-partition} algorithm since it allocates before
partitioning the balls.    First, allocate  $D$ balls randomly into
bins $\{\Ll\}\cup V$ with the restriction that $\Ll$ contains
exactly $\ell$ balls and each bin in $V$ contains at least $k$
balls. Then colour the balls u.a.r.\ with the following constraints:
\begin{description}
\item{(i)} exactly $im_i$ balls are coloured with $c_i$;

\item{(ii)} for each $i=2,\ldots,h$, the number of balls with colour
$c_i$ contained in $\Ll$ is exactly $l_i$.
\end{description}
Finally,  take u.a.r.\ a partition of the balls such that for each
$i=2,\ldots,h$, all balls with colour $c_i$ are partitioned into
parts of size $i$.

 To prove Theorem~\ref{t:CoreDensity}, we will convert the
problem to a question about $\mathcal{P}(V,{\bf M},{\bf L},k+1)$, in
particular the $(w,k+1)$-core of the hypergraph induced in the
obvious way by the bins containing at least $k+1$ balls.

A deletion algorithm   producing the $k$-core of a random multigraph
was analysed in ~\cite{CW}.
 The differential equation method~\cite{W2} was used to
analyse the size and the number of hyperedges of the final $k$-core.
The degree distribution of the $k$-core was shown to be a truncated
multinomial. We now extend this deletion algorithm to find the
$(w,k+1)$-core of $H$ in $\mathcal{M}_{\bar n,\bar m,h}$ and
$\mathcal{G}_{\bar n,\bar m,h}$. We describe the algorithm in the
setting of representing multihypergraphs using bins for vertices,
where each hyperedge   $x$ is a set $h(x)$ of $|x|$ balls. Initially
let $LV$ be the set of all light vertices/bins, and let $\overline
{LV}=V(H)\setminus LV$ be the set of heavy vertices. A {\em light}
ball is any ball contained in $LV$.
\begin{tabbing}
aaa\= aaa\= aaa\= \kill
{\bf RanCore Algorithm to obtain the $(w,k+1)$-core}\\
Input: an $h$-hypergraph $\orH$. Set $t:=0$.\\
\> While neither $LV$ nor $\overline {LV}$ is empty,\\
   \> \> $t:=t+1$;\\
   \> \>Remove all empty bins;\\
   \> \> U.a.r.\ choose a light ball $u$. Let $x$ be the hyperedge that contains $u$ and let $v$ be\\
   \> \>  the  vertex that contains $u$; \\
   \> \>If $|x|\ge h-w+2$, update $x$ with $x\setminus\{u\}$, \\
   \> \> \> otherwise,  remove this hyperedge $x$ from the current hypergraph. If any vertex \\
   \> \> \> $v'\in\overline{LV}$ becomes light, move $v'$ to $LV$ together with all balls in it;\\
\>If $LV$ is empty, ouput the remaining hypergraph, otherwise,
output the empty graph.
 \end{tabbing}

We will prove Proposition~\ref{uniformCore} and
Theorem~\ref{t:CoreDensity} by analysing the RanCore algorithm using
the partition-allocation model and the allocation-partition
algorithm which generates $\mathcal{P}(V,{\bf M},{\bf 0},k+1)$.
 Define
\begin{equation}
f_k(\mu)=\sum_{i\ge
k}e^{-\mu}\cdot\frac{\mu^i}{i!}=1-\sum_{i=0}^{k-1}e^{-\mu}\cdot\frac{\mu^i}{i!},\lab{eq:f}
\end{equation}
for any integer $k\ge 0$. By convention, define $f_{k}(\mu)=1$ for
any $k<0$. Let $Z_{(\ge k)}$ be a truncated Poisson random variable
with parameter $\lambda$ defined as follows.
\begin{equation}
\pr(Z_{(\ge
k)}=j)=\frac{e^{-\lambda}}{f_{k}(\lambda)}\cdot\frac{\lambda^j}{j!},\
\ \ \mbox{for any}\ j\ge k.\lab{eq:truncatedPoisson}
\end{equation}
Note that it follows that $\pr(Z_{(\ge k)}=j)=0$ whenever $j<k$. The
following proposition will be used in the proof of
Theorem~\ref{t:CoreDensity} and in Section~\ref{probabilistic}.
\begin{prop}\lab{p:lambda-mu}
  For
$\mu\ge k+2$, there exists a unique real $\lambda$ satisfying
$\lambda f_k(\lambda)=\mu f_{k+1}(\lambda)$.  Moreover,
$\mu\ge\lambda$, and if $\mu\ge ck$ for a fixed $c>1$, then
$\mu-\lambda\rightarrow 0$ and $f_k(\lambda)\rightarrow 1$ as
$k\to\infty$.
\end{prop}
\proof  Since $yf_k(y)/f_{k+1}(y)$ is monotonic in the domain $y>0$, as
shown in~\cite[Lemma 1]{PW}, there exist a unique $\lambda>0$ that satisfies
$
\lambda f_k(\lambda)=\mu f_{k+1}
$
as long as $\mu\ge \inf_{y>0}\{yf_k(y)/f_{k+1}(y)\}$. Clearly $f_k(1)/f_{k+1}(1)<k+2$, so $\mu\ge k+2$ suffices.
 Since
$f_k(\lambda)\ge f_{k+1}(\lambda)$ for all $k\ge -1$ by the
definition of the function $f_k(x)$ in~\eqn{eq:f}, it follows
directly that $\mu\ge \lambda$. For $k\to\infty$, we use well known
simple bounds on tails of the Poisson distribution. Set
$c=1+3\alpha$, where $\alpha>0$. If $\la<k+\alpha k$, then
\begin{eqnarray*}
\la f_k(\la)&\sim&\sum_{j=k}^{\lfloor k+2\alpha k \rfloor}
e^{-\la}\frac{\la^{j+1}}{j!} < (1+o(1))( k+2\alpha
k)\sum_{j=k}^{\lfloor k+2\alpha k \rfloor}
e^{-\la}\frac{\la^{j+1}}{(j+1)!}\\
  &<& ckf_{k+1}(\lambda)<\mu f_{k+1}(\lambda),
\end{eqnarray*}
a contradiction. So $\la\ge k+\alpha k$, whence $f_k(\lambda)$ and
$f_{k+1}(\lambda)$ are both $1+o(1/\la)$ and so $\la-\mu=o(1)$. \qed\ss

The following is essentially~\cite[Lemma 4.2]{W5}.
\begin{lemma}\lab{l:super}
Let $c>0$, $\delta$ be constants. Let $(Y_t)_{t\ge 1}$ be
independent random variables such that $|Y_t|\le c$ always and $\ex
Y_t\le \delta$ for all $t\ge 1$. Let $X_0=0$ and $X_t=\sum_{i\le t}
Y_i$ for all $t\ge 1$. Then for any $\eps>0$, a.a.s.\ $X_n \le
\delta n+\eps|\delta|n$. More precisely, $\pr(X_n\ge \delta
n+\eps|\delta| n)\le\exp(-\Omega(\eps^2 n))$.
\end{lemma}
\no {\bf Proof of Proposition~\ref{uniformCore}.\ } Consider an
element $\orH\in \mathcal{M}_{\bar n,{\bf m}}$ arising from $P\in
\mathcal{P}_{\bar n,{\bf m}}$, where ${\bf m} = (0,\ldots,0,\bar m)$
with $\bar m$ corresponding to the value of the coordinate $m_h$. If
we merge all the bins of $P$ containing   $k$ or fewer balls  into
one bin $\Ll$, we obtain in an obvious way an element $P'\in
\mathcal{P}(V,{\bf M},{\bf L},k+1)$ for an appropriate sequence
${\bf L}=(0,0,\ldots , L_0)$ where $L_0$ is the total degree of
light vertices (vertices with degree at most $k$) in $\orH$. Given
the parameters $(V,{\bf M},{\bf L},k+1)$, the number of $P$ that
will produce $P'$ is independent of $P'$. It follows that,
conditional on the total degree of the light vertices in $H$, this
generates $\mathcal{P}(V,{\bf M},{\bf L},k+1)$ with the correct
distribution. Moreover, the hypergraph induced   by the vertices of
$H$ of degree at least $k+1$ is also induced in the obvious way by
the heavy bins of $P'$. Hence, it suffices to study the
$(w,k+1)$-core of this hypergraph, conditional upon any feasible
${\bf L}$.  Because of the correspondence between
$\mathcal{P}(V,{\bf M},{\bf L},k+1)$ and the random
multihypergraphs, we sometimes call bins in $V$ vertices and the {\em
degree sequence} of $V$ denotes the sequence of numbers of balls
in bins in $V$.

To this end, we adapt the RanCore algorithm in the obvious way  to
be run on
  $\mathcal{P}(V,{\bf M},{\bf L},k+1)$, as follows.
  In each step $t$, the algorithm removes a ball, denoted
by $u$, u.a.r.\ chosen from all balls in $\Ll$.
  If the colour of $u$ is $c_{h-j}$ for
 $i<w-1$, the algorithm recolours the balls in the same part as $u$
 with the new colour $c_{h-j-1}$. If the colour of $u$ is
 $c_{h-w+1}$, the algorithm removes all balls contained in the same
 part as $u$, and if any heavy bin becomes light (i.e.\ the number of balls
 contained in it becomes at most $k$) because of the
 removal of balls, the bin is removed and the balls remaining in it are put into
 $\Ll$. This clearly treats the heavy bins of $P'$ in a corresponding way to
 RanCore treating the heavy vertices of $H$. Thus, the modified RanCore stops
 with a final partition-allocation that corresponds to the $(w,k+1)$-core of $H$, and this is what we will analyse.

For easier reference, let $g_t$ denote the random
partition-allocation derived after $t$ steps of this process.
 Let $V_t$ denote its set of heavy bins and let ${\bf
M}_t$ denote the class of sets $\{M_{t,h-w+1},\ldots,M_{t,h}\}$ such
that $M_{t,h-j}$ denotes the  set of partitioned balls with colour
$c_{h-j}$ in $g_t$ for $0\le j\le w-1$. Let
$m_{t,h-j}=|M_{t,h-j}|/(h-j)$ and ${\bf
m}=(m_{t,h-w+1},\ldots,m_{t,h})$. Let $L_{t,h-j}$ denote the number
of balls with colour $c_{h-j}$ in $\Ll$ and let ${\bf
L}=(L_{t,h-w+1},\ldots,L_{t,h})$. Let
$L_t=\sum_{j=0}^{w-1}L_{t,h-j}$. Initially, $g_0=P'$, $V_0=V$ etc.

There is a straightforward way to see,  by induction on $t$, that
the partition-allocation $g_t$, conditional on $V_t$, ${\bf M}_t$
and ${\bf L}_t$, is distributed as $\mathcal{P}(V_t,{\bf M}_t,{\bf
L}_t,k+1)$. We have already noted that this is true for $t=0$.
 For the inductive step,
 it suffices to note that that for any $t\ge 0$, conditional on the values of
 $V_t, {\bf M}_t, {\bf L}_t$ for $g_t$, the probability that $g_{t+1}$ is
 any particular member $g'$ of $V_{t+1}, {\bf M}_{t+1}, {\bf L}_{t+1}$
 does not depend on $g'$. This is because of three facts. Firstly, $g_t$
 is uniform conditional on the parameters at step $t$. Secondly,   the change
 in the parameters determines which type of step the algorithm is taking
 (e.g.\ if a heavy bin becomes light). Thirdly, each possibility for $g_{t+1}$
 is reachable from the same number of $g_t$ and, given the type of step occurring,
 each such transition has the same probability of occurring as step $t+1$.

Furthermore, it is easy to see that the degree distribution of the
$(w,k+1)$-core, conditional on the number of hyperedges of each
size, is truncated multinomial. This is because, for  any $V$ and
${\bf M}$, the allocation-partition algorithm which generates
$\mathcal{P}(V,{\bf M},{\bf 0},k+1)$ produces a truncated
multinomial distribution for the degrees of vertices in $V$. \qed\ss

The proof of Theorem~\ref{t:CoreDensity} uses the differential
equation method (d.e.\ method). In particular, we use the following special case
of~\cite[Theorem 6.1]{W5}. For each $n>0$
let a sequence of random vectors $(Y^{(1)}_t,\ldots,Y^{(l)}_t)_{0\le
t\le m}$ be defined on a probability space $\Omega_n$. (We suppress
the notation $n$.)  Let $U_t$ denote the history of the process up
to step $t$.

\begin{thm}\lab{t:de}
 Suppose that
there exists $C>0$ such that for each $i$, $Y^{(i)}<Cn$ always. Let
$\widehat \D\subset \mathbb{R}^{l+1}$ and let the stopping time $T$
be the minimum $t$ such that
$(t/n,Y_t^{(1)}/n,\ldots,Y_t^{(l)})\notin\widehat \D$. Assume
further that the following three hypotheses are satisfied.
\begin{description}
\item{(a)} (Boundedness hypothesis.) There exists a constant $C'>0$ such that for all $0\le t\le \min\{m,T\}$,
$|Y_{t+1}-Y_t|<C'$ always;

\item{(b)} (Trend hypothesis.) There exists functions $f_i$ for all $1\le i\le l$ such that for all $0\le t\le \min\{m,T\}$ and all $1\le i\le l$,
$$
\ex(Y_{t+1}^{(i)}-Y_t^{(i)}\mid
U_t)=f_i(t/n,Y_t^{(1)}/n,\ldots,Y_t^{(l)}/n)+o(1);
$$

\item{(c)} (Lipschitz hypothesis.) For every $1\le i\le l$, the functions $f_i$
are Lipschitz continuous in all their variables on a bounded
connected open set $\D$ where $\D$ contains the intersection of
$(t,z^{(1)},\ldots, z^{(l)}:t\ge 0)$ with some neighbourhood of
$(0,z^{(1)},\ldots, z^{(l)}:\pr(Y_0^{(i)}=z^{(l)}n,1\le i\le l)\neq
0\ \mbox{for some n})$.
\end{description}

Then the following conclusions hold.

\begin{description}
\item{(a)} For any $(0,\hat z^{(1)},\ldots,\hat z^{(l)})\in \D$, the
differential equation system
$$
\frac{d~z_i}{d~s}=f_i(s,z_1,\ldots,z_l), \ \ i=1,\ldots, l
$$
has a unique solution in $\D$ for $z_l:I\!\!R\to I\!\!R$ with the
initial conditions
$$
z_i(0)=\hat z^{(i)},\ \ i=1,\ldots,l,
$$
where the solution is extended arbitrarily close to the boundary of
$\D$.

\item{(b)} A.a.s.\
$$
Y_t^{(i)}=nz_i(t/n)+o(n)
$$
uniformly for all $0\le t\le
\min\{\sigma n, T\}$, where $\sigma$ is the supremum of all $x$ such
that the solution $(z^{(i)}(x))_{1\le i\le l}$ to the differential
equation system lies inside the domain $\D$.

\end{description}

\end{thm}
Our usage of a.a.s.\ in conjunction with other asymptotic notation such as $o()$ conforms to the conventions in~\cite{GTH}. For more details of the
method and proofs, readers can refer to~\cite[Theorem
1]{W2},~\cite[Theorem 5.1]{W5} and~\cite[Theorem 6.1]{W5}.  In our
case, $\Omega_n$ is the probability space of sequences of random
partition-allocations generated by running the RanCore algorithm on
graphs with $n$ vertices, where $t$ refers to the $t$-th step of the
algorithm and the $Y_t$  are variables defined during the
algorithm.)

The idea of the proof  of Theorem~\ref{t:CoreDensity} is, roughly
speaking, as follows.   We use the d.e.\ method to analyse the
asymptotic values of random variables defined on the random process
generated by the RanCore algorithm. The difficulty  arises from the
fact that the natural functions $f_i$ for our application are not
Lipschitz continuous at $x=0$. To avoid this, we artificially modify
the $f_i$ in a neighbourhood of the problem point, and show that the
solution to the new differential equation system coincides with the
original inside a domain $D_0$ which contains all points relevant to
the random process. Theorem~\ref{t:de} then applies to show that the
asymptotic values of the random variables inside $D_0$ are
approximated by the solution of the system, we analyse the random
variables when they leave  $D_0$. We show that provided $k$ is
sufficiently large, the algorithm then terminates quickly, which
allows us to estimate the size and density of $\coreH$.\ss

\no {\bf Proof of Theorem~\ref{t:CoreDensity}.\ } It was shown in
the proof of Proposition~\ref{uniformCore} that for every $t$,
conditional on the values of $V_t$, ${\bf M}_t$ and ${\bf L}_t$, the
partition-allocation $g_t$ is distributed as $\mathcal{P}(V_t,{\bf
M}_t,{\bf L}_t,k+1)$. \remove{Now we analyse the size of the
$(w,k+1)$-core.} After
 step $t$ of the RanCore  algorithm, define (or recall)  the following random variables:
\begin{eqnarray*}
B_t && \mbox{total number of balls remaining}\\
B_{t,h-j} && \mbox{number of balls coloured $c_{h-j}$}\\
A_{t,i} && \mbox{number of bins containing exactly $i$ balls}\\
 A_t  && \mbox{ $A_{t,k+1}$ (number of bins containing exactly $k+1$ balls)}\\
L_t && \mbox{number of light balls}\\
L_{t,h-j} && \mbox{number of light balls that are coloured
$c_{h-j}$}\\
H_{t,h-j} && \mbox{number of balls contained in heavy bins that
are coloured $c_{h-j}$}\\
 HV_t&& |V_t|, \mbox{ number of  heavy bins}
\end{eqnarray*}
and note that $H_{t,h-j} =B_{t,h-j}-L_{t,h-j}$.

Recall that $n$
and $m_{h-j}$ denote the number of vertices and the number of
hyperedges of size $h-j$ in $\coreH$, the $(w,k+1)$-core of $\orH$,
and $\bar\mu$ denotes its average degree.
 We show that if $\bar\mu\ge ck$ for some $c>1$ and
$k$ is sufficiently large, then a.a.s.\ $n\sim \alpha \bar n$ and
$m_{h-j}\sim \beta_{h-j}\bar n$ for some constants $\alpha>0$,
$\beta_{h-j}>0$ which are determined by the solution of the
  differential equation system given below, on a domain $D_0$ defined
below~\eqn{eq:DefSystem}.  In particular, we will show that
\begin{equation}
\alpha=z_{HV}(x^*), \ \
\beta_{h-j}=z_{H,h-j}(x^*)/(h-j),\lab{alphabeta}
\end{equation}
 where $x^*$ is the smallest positive root of
$z_L(x)=0$.

The d.e.\ method relies on a relation between solutions of a differential equation system and the random variables of the process under consideration. We will use subscripts of the real valued
functions to indicate their corresponding random variables.
For instance, the real function $z_{L,h-j}(x)$ is associated with
the random variable $L_{t,h-j}$. The differential equation system is as follows.
\begin{eqnarray}
z'_{L,h-j}(x)&=&\frac{z_{L,h-j}}{z_L}\left(-1-\frac{(h-j-1)z_{L,h-j}}{z_{B,h-j}}\right)\nonumber\\
&&+\frac{z_{L,h-w+1}}{z_L}\left(\frac{(h-w)z_{H,h-w+1}}{z_{B,h-w+1}}\cdot\frac{(k+1)z_A}{z_B-z_L}\cdot
k\cdot\frac{z_{H,h-j}}{z_B-z_L}\right)\nonumber\\
&&+\frac{z_{L,h-j+1}}{z_L}\frac{(h-j)z_{L,h-j+1}}{z_{B,h-j+1}},
\ \ j=1,\ldots,w-1,\lab{eq:DL}\\
z'_{H,h-j}(x)&=&\frac{z_{L,h-j}}{z_L}\left(-\frac{(h-j-1)z_{H,h-j}}{z_{B,h-j}}\right)\nonumber\\
&&-\frac{z_{L,h-w+1}}{z_L}\left(\frac{(h-w)z_{H,h-w+1}}{z_{B,h-w+1}}\cdot\frac{(k+1)z_A}{z_B-z_L}\cdot
k\cdot\frac{z_{H,h-j}}{z_B-z_L}\right)\nonumber\\
&&+\frac{z_{L,h-j+1}}{z_L}\frac{(h-j)z_{H,h-j+1}}{z_{B,h-j+1}}, \ \ j=1,\ldots,w-1,\lab{eq:DH}\\
z_L'(x)&=&-1+\frac{z_{L,h-w+1}}{z_L}\bigg(-\frac{(h-w)z_{L,h-w+1}}{z_{B,h-w+1}}+(h-w)k\cdot\frac{z_{H,h-w+1}}{z_{B,h-w+1}}\cdot\frac{(k+1)z_A}{z_B-z_L}\bigg)\lab{eq:DD}\\
z_B'(x)&=&-1-\frac{(h-w)z_{L,h-w+1}}{z_L}\lab{eq:DB}\\
z_{HV}'(x)&=&-\frac{z_{L,h-w+1}}{z_L}\frac{(h-w)z_{H,h-w+1}}{z_{B,h-w+1}}\cdot\frac{(k+1)z_A}{z_B-z_L}\lab{eq:DHV}\\
\lambda'(x)&=&\frac{((z_B'-z_L')z_{HV}-(z_B-z_L)z_{HV}')f_{k+1}(\lambda)}{z_{HV}^2(f_k(\lambda)+\lambda
e^{-\lambda}\cdot
\frac{\lambda^{k-1}}{(k-1)!}-\frac{z_B-z_L}{z_{HV}}\cdot
e^{-\lambda}\cdot \frac{\lambda^k}{k!})}\lab{eq:Dmu}\\
z_{L,h}(x)&=&z_L(x)-\sum_{i=1}^{w-1}z_{L,h-j}(x),\ \ \ z_{H,h}(x)=z_B(x)-z_L(x)-\sum_{i=1}^{w-1}z_{H,h-j}(x),\lab{eq:LH}\\
z_{B,h-j}(x)&=&z_{L,h-j}(x)+z_{H,h-j}(x), \ \ \mbox{for every}\ 0\le j\le w-1,\lab{eq:Bi}\\
z_A(x)&=&\frac{\lambda(x)^{k+1}}{e^{\lambda(x)}(k+1)!f_{k+1}(\lambda(x))}z_{HV}(x),\lab{eq:DA}
\end{eqnarray}
where $f_k(\lambda)$ was defined in~\eqn{eq:f}. The initial
conditions are
\begin{eqnarray}
z_B(0)&=&\bar\mu,\   z_{L,h-j}(0)=0, \   z_{H,h-j}(0)=0, \ \mbox{for
all}\ 1\le j\le w-1, \lab{eq:initial} \\
z_L(0)&=&\bar\mu(1-f_k(\bar\mu)), \
z_{HV}(0)=1-\exp(-\bar\mu)\sum_{i=0}^k \bar\mu^i/i!,\
\lambda(0)=\bar\mu. \lab{eq:DefSystem}
\end{eqnarray}
Let $D_0$ be the domain which contains all points such that
$x\in\mathbb{R}$, $0\le z_{L,h-j}\le z_{B,h-j}$, $z_{L,h-j}\le z_L$
for all $0\le j\le w-1$, $z_L>0$, $z_B-z_L>0$,   $z_{HV}>0$
and $(z_B-z_L)/z_{HV}> k+2$.   We will call the right hand sides
of~(\ref{eq:DL}--\ref{eq:Dmu}) the {\em
derivative  functions}, and at present we   regard them to be only defined in
$D_0$.  It is straightforward to check that for any point ${\bf z^*}\in D_0$ such that $z_{L,h-j}=z_{B,h-j}=0$, the functions specified
in the right hand sides of~(\ref{eq:DL}--\ref{eq:DHV}) tend to $0$
when $\bf z$ approaches
  $\bf z^*$ from the interior of  $D_0$.  For example, note that the term
$$
\frac{z_{L,h-j}}{z_L}\cdot(h-j-1)\cdot\frac{z_{L,h-j}}{z_{B,h-j}}
$$
on the right hand side of~\eqn{eq:DL}
is bounded above by $(h-j-1) z_{L,h-j}/{z_L}$
since $|z_{L,h-j}/z_{B,h-j}|\le 1$ when ${\bf z}\in D_0$. Hence, it tends to 0 if ${\bf z} \to {\bf z}^*$. The same applies to similar terms in~\eqn{eq:DL}--\eqn{eq:DHV}.
  As part of our definition of the differential equation system~\eqn{eq:DL}--\eqn{eq:DefSystem}, we now declare the values of
these terms at such points $\bf z^*$ to be $0$.

In applying Theorem~\ref{t:de}, the variable $x$ will be
associated with $t/\bar n$. As mentioned above, the variable $z_{L,h-j}(x)$
is associated with the variable  $L_{t,h-j}/\bar n$, which we call
the {\em scaled version} of the random variable  $L_{t,h-j}$. We do
the same for the other random variables, and call $t/\bar n$ the
scaled version of $t$.

There are two kinds of problems with the Lipschitz property required
in Theorem~\ref{t:de} (c). The first is caused by terms in the
equations with denominators $z_{L}$ or  $z_B-z_L$ appearing in the
derivative functions, which are potentially 0, causing
singularities. These are relatively easy to take care of since they
do not become small until near the end of the process. For any fixed
constant $\eps>0$, define $\hD(\eps)$ to be the connected subset of
$D_0$ obtained by restricting to ${\bf z}$ such that $z_{L}>\eps$
and $z_B-z_L>\eps$. We will basically restrict consideration to
points in $\hD(\eps)$.  Let $T$ be the (stopping) time that the
vector of scaled random variables leaves $\hD(\eps)$. Let $t\wedge
T$ denote $\min\{t,T\}$. The conclusion of Theorem~\ref{t:de} will
give information on the scaled random variables up to the step when
they reach the boundary of the domain $\hD(\eps)$. This gives us
information about $(g_{t\wedge T})_{0\le t\le \tau}$, where $g_t$ is
the partition-allocation obtained after step $t$. At that point we
will need some further observations to show that the process
finishes soon afterwards.

The second type of problem comes from denominators containing $z_{B,h-j}$, which can be 0 even right at the start of the process. This poses a difficulty since the theorem requires the derivative functions to be Lipschitz in an open domain containing the starting point. To deal with this, we will, at an appropriate point below, extend the
differential equations into a larger connected open domain $D\supset
D_0$, and correspondingly extend $\hD(\eps)$ to $D(\eps)$. We will
actually
apply Theorem~\ref{t:de} with $\D=D(\eps)$ and
$\widehat\D=\hD(\eps)$.

We first verify hypotheses (a) and (b) of Theorem~\ref{t:de}, which are unrelated to the choice of $D(\eps)$.
It is easy to see that the change of each random variable in every step of the algorithm is
bounded.  This is because in every step, the number of balls
deleted (or recoloured, or moved from heavy bins to the light bin $\Ll$) is
bounded.
Thus, Theorem~\ref{t:de}(a) clearly holds.

To
verify hypothesis (b), we will need to show that the expected
one-step change of each random variable, such as $L_{t,h-j}$, can be
approximated to within $o(1)$ error by some function of the scaled
variables. Replacing the scaled variables in these functions  by
their associated real variables will give the derivative functions
in~\eqn{eq:DL}--\eqn{eq:DHV}.

 Let $g_t$ be the partition
allocation obtained after step $t$. At step $t+1$, a
partition-allocation $g_{t+1}$ is to be obtained by applying the
RanCore algorithm to $g_t$. Let $v$ be the ball randomly chosen by
the algorithm from $\Ll$. Let $C(v)$ be the colour of $v$, so
$C(v)=h-j$ for some $j$. If $j<w-1$, the algorithm removes another
$h-j-1$ balls that are uniformly distributed among all balls with
colour $c_{h-j}$ since $g_t\in \mathcal{P}(V_t,{\bf M}_t,{\bf
L}_t,k+1)$ as proved in Proposition~\ref{uniformCore}.
 If $j=w-1$, then the algorithm removes
 $v$ together with another $h-w$ balls which are chosen u.a.r.\ from all
balls of colour $c_{h-w+1}$. If the removal of the $h-w$ balls
results in some heavy bins turning into light bins, these bins are
removed and the balls remaining in these bins are put into $\Ll$.

 Now we estimate the expected value of
$L_{t+1,h-j}-L_{t,h-j}$ for any $1\le j\le w-1$ and for any $0\le
t<\tau$ conditional on $V_{t}$, ${\bf M}_{t}$, ${\bf L}_{t}$ and the
event $g_{t}\in\mathcal{P}(V_{t},{\bf M}_{t},{\bf L}_{t},k+1)$. Given $j$, the
probability that
 $C(v)=c_{h-j}$ is $L_{t,h-j}/L_t$.  If $C(v)=c_{h-j}$, one ball
 of colour $c_{h-j}$ contained in $\Ll$ is removed, and
 another $h-j-1$ balls of colour $c_{h-j}$ are recoloured with
 $c_{h-j-1}$ (or removed if $j=w-1$). So the expected number of those
 balls
 that are contained in $\Ll$ is
 $$
\frac{(h-j-1)L_{t,h-j}}{B_{t,h-j}}(1+o(1)),
 $$
 provided $B_{t,h-j}\ge \log n$ (say).
 Hence
$$
\frac{L_{t,h-j}}{L_t}\left(-1-\frac{(h-j-1)L_{t,h-j}}{B_{t,h-j}}\right)+o(1)
$$
is the negative contribution to $\ex(L_{t+1,h-j}-L_{t,h-j}\mid
V_{t},{\bf M}_{t},{\bf L}_{t}, g_{t}\in\mathcal{P}(V_{t},{\bf
M}_{t},{\bf L}_{t},k+1))$. Note that we reach the same conclusion if
$B_{t,h-j}<\log n$ because in that case
$$
L_{t,h-j}/L_t\le B_{t,h-j}/L_t<\log n/\eps n=o(1).
$$
The positive contribution to $\ex(L_{t+1,h-j}-L_{t,h-j}\mid
V_{t},{\bf M}_{t},{\bf L}_{t}, g_{t}\in\mathcal{P}(V_{t},{\bf
M}_{t},{\bf L}_{t},k+1))$ comes from the following two cases.

\no {\em Case 1: } $C(v)=c_{h-w+1}$. Here,  the algorithm removes
$v$ and another $h-w$ balls of colour $c_{h-w+1}$.
$\pr(C(v)=h-w+1)=L_{t,h-w+1}/L_t$.  We first note  that, for any
$2\le i\le h-w$, the contribution from the case that   $i$ of the
$h-w$ removed balls lie in a bin containing at most $k+i$ balls is
at most $\big((k+i)/(B_t-L_t)\big)^{i-1} =o(1),
 $
since the  definition of $\widehat\D=\hD(\eps)$  ensures that the denominator is at least $\eps n$ for $t\le T$.

 It only remains to consider the contribution
from the case that a ball in a bin containing exactly $k+1$ balls is
removed. For each ball removed, the probability that it is in a bin
containing exactly $k+1$ balls is
$$
\frac{H_{t,h-w+1}}{B_{t,h-w+1}}\cdot\frac{(k+1)A_{t,k+1}}{B_t-L_t}+o(1).
$$
 The removal of such a ball causes the bin to become light.
 Since $g_{t}\in\mathcal{P}(V_{t},{\bf
M}_{t},{\bf L}_{t},k+1)$, the balls of each colour are uniformly
distributed among all balls in the heavy bins, and thus the expected
number of balls of colour $c_{h-j}$, for $0\le j\le w-1$, among the
remaining $k$ balls in the bin  is
$$
k\cdot\frac{H_{t,h-j}}{B_t-L_t}+o(1).
$$
In total, $h-w$ balls of colour $c_{h-w+1}$ are removed, other than
$v$. Hence the expected contribution to
$\ex(L_{t+1,h-j}-L_{t,h-j}\mid V_{t},{\bf M}_{t},{\bf L}_{t},
g_{t}\in\mathcal{P}(V_{t},{\bf M}_{t},{\bf L}_{t},k+1))$ is
$$
(h-w)\cdot\frac{L_{t,h-w+1}}{L_t}\cdot\frac{H_{t,h-w+1}}{B_{t,h-w+1}}\cdot\frac{(k+1)A_{t,k+1}}{B_t-L_t}\cdot
k\cdot\frac{H_{t,h-j}}{B_t-L_t}+o(1).
$$

\no {\em Case 2: } $C(v)=c_{h-j+1}$. The algorithm removes $v$, chooses another $h-j$ balls u.a.r.\ from those of colour
$c_{h-j+1}$, and recolours them with $c_{h-j}$. Since
$\pr(C(v)=c_{h-j+1})=L_{t,h-j+1}/L_t$,  conditional on
$C(v)=c_{h-j+1}$, the expected number of balls of colour $c_{h-j+1}$
that are in the light bins and are recoloured   is
$$
(h-j)\cdot\frac{L_{t,h-j+1}}{B_{t,h-j+1}}+o(1),
$$
provided $B_{t,h-j+1}\ge \log n$. Hence the positive contribution to
$\ex(L_{t+1,h-j}-L_{t,h-j}\mid V_{t},{\bf M}_{t},{\bf L}_{t},
g_{t}\in\mathcal{P}(V_{t},{\bf M}_{t},{\bf L}_{t},k+1))$ is
$$
\frac{L_{t,h-j+1}}{L_t}\cdot(h-j)\cdot\frac{L_{t,h-j+1}}{B_{t,h-j+1}}+o(1)
$$
in this case. The same conclusion holds when $B_{t,h-j+1}<\log n$
for the same reason as discussed before. Therefore
\begin{eqnarray}
&&\ex(L_{t+1,h-j}-L_{t,h-j}\mid V_{t},{\bf M}_{t},{\bf L}_{t},
g_{t}\in\mathcal{P}(V_{t},{\bf
M}_{t},{\bf L}_{t},k+1))\nonumber\\
&&\hspace{.3cm}=\frac{L_{t,h-j}}{L_t}\left(-1-\frac{(h-j-1)L_{t,h-j}}{B_{t,h-j}}\right)+\frac{L_{t,h-j+1}}{L_t}\cdot
\frac{(h-j)L_{t,h-j+1}}{B_{t,h-j+1}}\nonumber\\
&&\hspace{.7cm}+\frac{L_{t,h-w+1}}{L_t}\left(\frac{(h-w)H_{t,h-w+1}}{B_{t,h-w+1}}\cdot\frac{(k+1)A_{t,k+1}}{B_t-L_t}\cdot
k\cdot\frac{H_{t,h-j}}{B_t-L_t}\right)+o(1),\lab{eq:BL}
\end{eqnarray}
for $j=1,\ldots,w-1$. Replacing the random variables in the right hand side of~\eqn{eq:BL} by their associated real variables (noting that the scaling cancels out) gives
the right hand side of~\eqn{eq:DL}.
Using a similar approach to computing the expected changes of $H_{t,h-j}$, $B_{t}$, $D_{t}$, $HV_{t}$,
 conditional on   $V_{t}$, ${\bf M}_{t}$, ${\bf L}_{t}$ and the event $g_{t}\in\mathcal{P}(V_{t},{\bf
M}_{t},{\bf L}_{t},k+1)$,   we easily obtain the derivative
functions in~\eqn{eq:DH}--\eqn{eq:DHV}. The  equations
\begin{eqnarray*}
L_{t,h}&=&L_t-\sum_{i=1}^{w-1}L_{t,h-j},\ \ \
H_{t,h}=B_t-L_t-\sum_{i=1}^{w-1}H_{t,h-j},\\
B_{t,h-j}&=&L_{t,h-j}+H_{t,h-j},\ \ \ \mbox{for every}\ h-w+1\le
j\le h
\end{eqnarray*}
are obvious and lead to~\eqn{eq:LH} and~\eqn{eq:Bi}.

Let $\mu_t$ denote $(B_t-L_t)/HV_t$, the average degree of heavy
vertices after step $t$. Correspondingly we define a function $\mu(x)$
associated with the random variable $\mu_t$ to be
\begin{equation}
\mu(x)=(z_B(x)-z_L(x))/z_{HV}(x).\lab{mu}
\end{equation}
Then by Proposition~\ref{p:lambda-mu}, we may define $\la(x)$ by
\begin{equation}
\lambda(x)f_k(\lambda(x))=\mu(x)f_{k+1}(\lambda(x))\lab{lambdaofx}
\end{equation}
provided that $\mu(x)>k+2$, which is guaranteed inside $D_0(\eps)$. Let
$\lambda_t=\lambda(t/\bar n)$, so that  $\lambda_t$
is the unique positive root of
\begin{equation}
\frac{\lambda_tf_k(\lambda_t)}{f_{k+1}(\lambda_t)}-\mu_t=0.\lab{root}
\end{equation}

 Since $g_t\in\mathcal{P}(V_t,{\bf M}_t,{\bf L}_t,k+1)$ for every $t$,
 by considering the allocation-partition algorithm that generates
 $\mathcal{P}(V_t,{\bf M}_t,{\bf L}_t,k+1)$, the degree sequence of the
 heavy vertices has the truncated multinomial distribution. Hence, by~\cite[Lemma
1]{CW},
\begin{equation}
A_{t,k+1}\sim
\frac{e^{-\lambda_t}\lambda_t^{k+1}}{(k+1)!f_{k+1}(\lambda_t)}HV_t,\lab{eq:A}
\end{equation}
where  $\lambda_t$ satisfies~\eqn{root}. This
gives~\eqn{eq:DefSystem}.

Now~\eqn{eq:Dmu}, which gives the derivative of $\lambda(x)$, follows by
taking the derivative of both sides of~\eqn{lambdaofx},
$$
\lambda'(x) f_k(\lambda(x))+\lambda(x)
\frac{df_k(\lambda)}{d\lambda}\bigg|_{\lambda=\lambda(x)}
\lambda'(x)=\mu'(x) f_{k+1}(\lambda(x))+\mu(x)
\frac{df_{k+1}(\lambda)}{d\lambda}\bigg|_{\lambda=\lambda(x)}\lambda'(x),
$$
where,  by the definitions of   $f_k(\lambda)$ in~\eqn{eq:f} and
$\mu(x)$ in~\eqn{mu},
$$
\frac{df_k(\lambda)}{d\lambda}=e^{-\lambda}\frac{\lambda^{k-1}}{(k-1)!}\
\ \ \mbox{and}\ \ \
\mu'(x)=\frac{(z_B'-z_L')z_{HV}-(z_B-z_L)z_{HV}'}{z_{HV}^2}.
$$


Now we justify hypothesis (c).
We will first extend the derivative functions (which, up until this
point, we restricted to $D_0$) into a larger domain $D$, which
defines an extended d.e.\ system,  and show that these extended
functions are continuous and Lipschitz inside an open domain
$D(\eps)$ extended from $\hD(\eps)$. Later we will show  that the solution of
 the extended d.e.\
system for $0\le x\le T/\bar n$, with the same initial conditions as
the original system, is contained inside the domain $\hD(\eps)$ and
is thus the solution to the original d.e.\ system.

 We begin with the domain $\hD(\eps)$,
which was defined by restricting the points in $D_0$ to $z_{L}>\eps$
and $z_B-z_L>\eps$. \remove{******Let $f_{L,h-j}$ denote the
derivative function on the right hand side of~\eqn{eq:DL}; the other
ones are similar.****} Recalling our treatment of the possible
singularity $z_{L,h-j}=z_{B,h-j}=0$ just after~\eqn{eq:DefSystem},
each derivative function   is continuous in $\hD(\eps)$. The only
potential problems for the Lipschitz property  are the constant
multiples of the function
\begin{equation}
f(z_{L,h-j},z_{L},z_{B,h-j})=\frac{z_{L,h-j}}{z_{L}}\cdot\frac{z_{L,h-j}}{z_{B,h-j}}.\lab{singular}
\end{equation}
However, recalling that  $z_{L}>\eps$, $0\le z_{L,h-j}\le z_{B,h-j}$ and $z_{L,h-j}\le
z_L$ in $\hD(\eps)$, we have that the partial derivatives of
$f(z_{L,h-j},z_{L},z_{B,h-j})$
with respect to $z_{L,h-j}$, $z_{L}$ and $z_{B,h-j}$ are all $O(1/\eps)$,
\remove{ 
\begin{eqnarray}
&&\Big|\frac{\partial f}{\partial z_{L,h-j}}(z_{L,h-j},z_{L},z_{B,h-j})\Big|=\frac{2z_{L,h-j}}{z_{L}z_{B,h-j}}\le\frac{2}{z_L}<\frac{2}{\eps},\lab{part1}\\
&&\Big|\frac{\partial f}{\partial z_{L}}(z_{L,h-j},z_{L},z_{B,h-j})\Big|=\frac{z_{L,h-j}^2}{z_{L}^2z_{B,h-j}}<\frac{1}{\eps},\lab{part2}\\
&&\Big|\frac{\partial f}{\partial
z_{B,h-j}}(z_{L,h-j},z_{L},z_{B,h-j})\Big|=\frac{z_{L,h-j}^2}{z_{L}z_{B,h-j}^2}<\frac{1}{\eps}.\lab{part3}
\end{eqnarray}
}
from which it follows that the derivative functions are Lipschitz in
$\hD(\eps)$. \remove{****** \nickv{Changed and abbreviated. We have
now given a good proof sketch earlier so don't need to say so much
again. Also I think the order has been improved. It seems that  in
several places there are reminders of things that do need reminding
somewhere but are not so relevant at that particular place, which
interrupted the flow and resulted in more reminders being needed for
the other things! }********}

 Let ${\bf z}_0$ denote the initial
condition vector given by~\eqn{eq:initial} and~\eqn{eq:DefSystem}:
$x=0$, $z_L=\bar\mu(1-f_k(\bar\mu))$, $z_B=\bar\mu$,
$z_{HV}=1-\exp(-\bar\mu)\sum_{i=0}^k \bar\mu^i/i!$,
$z_{L,h-j}=z_{B,h-j}=0$ for all $1\le j\le w-1$. Note that    ${\bf
z}_0 $ lies on the boundary of both $D_0$ and $\hD(\eps)$.
Define
$D:=\{(x,z_{L,h-w+1},\ldots,z_{L,h-1},z_{B,h-w+1},\ldots,z_{B,h-1},z_L,z_B,z_{HV}):z_L>0,
z_B-z_L>0, z_B-z_L>(k+2)z_{HV}\}$, and let $D(\eps)$ be the domain obtained by restricting
points in $D$ to those with $z_L>\eps$ and $z_B-z_L>\eps$. Thus
$D(\eps)$ is the corresponding extension of $\hD(\eps)$.
 Clearly ${\bf z}_0$ is an interior point in $D$ and $D(\eps)$.
To extend the derivative functions to $D$, it is enough to extend the function
$f$ in~\eqn{singular}. Define
\begin{equation}
f^*(z_{L,h-j},z_{B,h-j},z_L)=\left\{\begin{array}{ll}
f(z_{L,h-j},z_{B,h-j},z_L) &\mbox{if}\ 0\le z_{L,h-j}\le z_{B,h-j},
z_{B,h-j}>0\\
0& \mbox{if}\ z_{L,h-j}=z_{B,h-j}=0,\\
z_{B,h-j}/z_L& \mbox{if}\ z_{L,h-j}>z_{B,h-j}\ge 0\\
f(|z_{L,h-j}|,|z_{B,h-j}|,z_L)&\mbox{otherwise}.
\end{array}\right.\lab{eq:extension}
\end{equation}
%
%
%
%
%
 We have already shown that $f$ is
Lipschitz continuous on $\hD(\eps)$, which is the first case
of~\eqn{eq:extension}. Since $z_L>\eps$ and $z_B-z_L>\eps$ in
$D(\eps)$, $f^*$ is Lipschitz continuous on $D(\eps)$. Hence, if we
modify the differential equation
system~\eqn{eq:DL}--\eqn{eq:DefSystem} by replacing each expression
equivalent to $f$ by $f^*$, we obtain derivative functions that are
Lipschitz continuous in the open domain $D(\eps)$. Thus, hypothesis
(c) holds for this system, which we call  the {\em extended
differential equation system}.

 We may now
apply Theorem~\ref{t:de}, to deduce that a.a.s.\ uniformly for every
$0\le t\le T$, $L_{t}=\bar nz_L(t/\bar n)+o(\bar n)$, and the same
applies to all the other random variables under consideration. We
claim that the stopping time $T$ coincides with the time at which
$L_t/\bar n$ or $(B_t-L_t)/\bar n$ decreases to $\eps$. This follows
by the following two observations, whose verifications are only sketched here since they require  straightforward
analysis. (See~\cite[pp.\ 86,87]{G3} for details.)

(i) The solution of the extended differential equation system is
interior to $\hD(\eps)$ for all sufficiently small $x> 0$.  For
instance, all functions taking the value $0$ at $x=0$ have positive
derivatives for sufficiently small $x>0$. Thus, these functions
become positive for any sufficiently small $x$ and thus the solution
is inside $\hD(\eps)$.)

(ii) Once the solution is interior to $\hD(\eps)$, the only
boundaries of $\hD(\eps)$ it can reach are  $z_{L}=\eps$,
$z_B-z_L=\eps$, $z_{HV}=0$ and $(z_B-z_L)/z_{HV}= k+2$. The other
boundaries of this domain are $z_{L,h-j}=0$, $z_{L,h-j}= z_{B,h-j}$
(i.e.\ $z_{H,h-j}=0$), and $z_{L,h-j} = z_L$ for any $j\ge 0$. For
example, it cannot reach
  $z_{L,h-j}=0$   because  the
only negative contribution to the derivative of $z_{L,h-j}$ is
proportional to $z_{L,h-j}$ itself. In view of this,  $z_{L,h-j} < z_L$ for any $0\le j\le w-1$ and the last-listed boundary cannot be reached.

Let $x(\eps)$  be the smallest value of $x$ such that $z_L(x)=\eps$,
$z_B(x)-z_L(x)=\eps$, $z_{HV}(x)=0$ or
$z_B(x)-z_L(x)=(k+2)z_{HV}(x)$, i.e., $\mu=k+2$, considering the definition~\eqn{mu}.  Then the solution to the extended
differential equation system  for all $0\le x\le x(\eps)$ is also
the solution to the original differential equation system. Let $x^*$
be the smallest real number such that $z_L(x^*)=0$,
$z_B(x^*)-z_L(x^*)=0$, $z_{HV}=0$ or $\mu=(k+2)$.
Then  the solution of the original differential equation system can
be extended arbitrarily close to $x^*$.

By the theorem's hypothesis, $\bar\mu\ge ck$ for some $c>1$. We next
show that for sufficiently large $k$ (depending on the value of
$c$), the function $z_L(x)$ reaches $0$ before $z_B(x)-z_L(x)$ or
$z_{HV}(x)$ reach $0$ or $\mu$ reaches $k+2$,
and we also provide an upper bound of the value of $x^*$. Let
$z_H(x)=z_B(x)-z_L(x)$. Clearly $z_L(x)=\sum_{i=0}^{w-1}
z_{L,h-j}(x)$ and $z_H(x)=\sum_{i=0}^{w-1} z_{H,h-j}(x)$.
So~\eqn{eq:DD},~\eqn{eq:DB} and~\eqn{eq:DHV} immediately lead to
\begin{eqnarray}
z_L'(x)\le -1+\frac{hk(k+1)z_A}{z_H},\ \   z_H'\ge
-h-\frac{hk(k+1)z_A}{z_H},\ \ z_{HV}'(x)\ge -\frac{h(k+1)z_A}{z_H}.
\lab{eq:derivative}
\end{eqnarray}
Let $\delta=(1-f_k(\bar\mu))\bar\mu$.  Then the initial conditions
give $z_L(0)=\delta$ and $z_H(0)=\bar\mu-\delta$. Since $\bar\mu\ge
ck$ for some $c>1$, $\delta=\exp(-\Omega_{c}(k))$.  By Proposition~\ref{p:lambda-mu}, we may assume that as
long as $\mu(x)\ge c'k$ for some $c'>1$ and $k$ sufficiently large,
$\lambda(x)$, is well defined by~\eqn{lambdaofx},   and $|\mu(x)-\lambda(x)|\le 1$, which
implies that $ z_A(x)/z_H(x)=\exp(-\Omega_{c'}(k))$ by~\eqn{eq:A}.
We next observe that $\mu(0)\ge\lambda(0)=\bar\mu$
by~\eqn{eq:DefSystem} and Proposition~\ref{p:lambda-mu}. Let
$[0,x_0]$ be an interval such that $\mu(x)\ge \bar\mu-4h$ for all
$0\le x\le x_0$. Certainly $\mu(x)\ge c'k$ for some $c'>1$ for all
$0\le x\le x_0$. We may choose $k$ sufficiently large (depending only
on the value of $c'$)  that $\delta\le 1$  and for all $0\le
x\le x_0$ we have $|\lambda(x)-\mu(x)|\le 1$,
$hk(k+1)z_A(x)/z_H(x)\le 1/2$ and $h(k+1)z_A(x)/z_H(x)\le 1/8$. Then
for all $0\le x<x_0$
\begin{equation}
z_L'(x)\le -1/2, \quad 0\ge z_H'(x)\ge -h-1/2, \quad z_{HV}'(x)\ge
-1/8. \lab{upper}
\end{equation}
Note that  $z_{HV}(0)\le 1$, and  $z'_{HV}(x)<0$ from~\eqn{eq:DHV}.
Thus $\mu(x)=z_H(x)/z_{HV}(x)\ge z_H(x)$ for any $0\le x<x^*$.
Hence, provided $z_H(x)\ge \bar\mu-4h$, we have $\mu(x)\ge
\bar\mu-4h$ and so the inequalities~\eqn{upper} hold. Then
$z_H(x)\ge z_{H}(0)+x\left(-h-\frac{1}{2}\right)=
\bar\mu-\delta+x(-h-\frac{1}{2})> \bar\mu-4h$ provided $x\le
3\delta$ say, since $\delta$ is arbitrarily small for large $k$. It
follows that $\mu(x)\ge z_H(x)\ge \bar\mu-4h$ for $x\le \min\{x^*,
3\delta\}$. Thus we may choose   $x_0\ge\min\{x^*, 3\delta\}$, and
so~\eqn{upper} implies, for any $0\le x<\min\{x^*, 3\delta\}$, that
\begin{equation}
z_{HV}(x)\ge z_{HV}(0)-\frac{3\delta}{8}>0,\ \ \ z_L(x)\le
\delta-\frac{x}{2}.\lab{eq:L}
\end{equation}
So $x^*<3\delta$ and $z_{HV}(x^*)>0$, since otherwise $3\delta\le
x^*$ and $z_L(3\delta)\le \delta-3\delta/2<0$, contradicting the
definition of $x^*$.  Combining this with $\mu(x)\ge z_H(x)\ge
\bar\mu-4h$, which is greater than $k+2$ for sufficiently large $k$,
we conclude that $z_L(x)$ reaches $0$ before $z_H(x)$ or $z_{HV}(x)$
reaches $0$ and before $\mu(x)$ reaches $k+2$ (in fact, before
$\mu(x)$ reaches $\bar\mu-4h$) and $x^*\le 3\delta$. We also have
that $z_L'(x)\le -1/2$ for all $x<x^*$.

  For notational convenience, define the following
limits from below (which we know to exist from the above bounds on
the functions and their derivatives):
\begin{equation}
z_{H,h-j}(x^*):=\lim_{x\to (x^*)^-}z_{H,h-j}(x)\ \ \mbox{and}\ \
z_{HV}(x^*):=\lim_{x\to (x^*)^-}z_{HV}(x).\lab{limit}
\end{equation}
Note that this definition yields continuous functions
$z_{H,h-j}(x)$ and $z_{HV}(x)$ on the closed interval $[0,x^*]$.

Given any sufficiently small $\eps>0$, let $x(\eps)$ be the root of
$z_L(x)=\eps$ and let $t(\eps)=\lfloor x(\eps)\bar n\rfloor$. Let
$Y_t$ denote  any of the random variables $H_{t,h-j}$ or $HV_{t}$,
 and $y(x)$ its associated real function.
We have shown that a.a.s.
\bel{conc}
 Y_{\lfloor x \bar n\rfloor }=\bar n
y (x )+o(\bar n)
\ee
 for $0\le x\le x(\eps)$.
Also, we have $|Y_{\lfloor x\bar n\rfloor}-Y_{\lfloor x(\eps)\bar
n\rfloor}|=O((x-x(\eps))\bar n)$ for all $x(\eps)\le x\le x^*$ since
the change of each variable in every step is bounded by $O(1)$. Let
$\delta_1(\eps)$ denote the number of light balls remaining at step
$t(\eps)$. Then $\delta_1(\eps)=\eps\bar n+o(\bar n)$. Applying
Lemma~\ref{l:super} with $X_0=L_{t(\eps)}$,
$X_n=L_{t(\eps)+4\delta_1(\eps)}$, $n=4\delta_1(\eps)$,
$\delta=-1/2$ and $c=h$, we have a.a.s.\
$L_{t(\eps)+4\delta_1(\eps)}\le
\delta_1(\eps)-(4\delta_1(\eps)/2)/2= 0$. Hence,  the time $\tau$
that the RanCore algorithm terminates a.a.s.\ satisfies $\tau \le
t(\eps)+4\delta_1(\eps)$. If it terminates before $\bar n x^*$, we
may artificially let it run to that point, with the variables
remaining static, thereby defining them on the interval $t\le \lfloor  x^*\bar n  \rceil$.
 Then, letting $\eps \to 0$ shows that the
conclusion~\eqn{conc} above applies for $0\le x\le x^*$,  noting that  the
function $y(x)$ is continuous on   $[0,x^*]$ as noted below~\eqn{limit}. We may also conclude, since
$\delta_1(\eps)\to 0$ as $\eps\to 0$, that \bel{tauvalue} \tau =
x^*\bar n +o(n)\ a.a.s. \ee

 In particular,
we conclude that a.a.s.\ $H_{\lfloor x\bar n\rfloor,h-j}=\bar n
z_{H,h-j}(x)+o(\bar n)$ and $HV_{\lfloor x\bar n\rfloor}=\bar n
z_{HV}(x^*)+o(\bar n)$. Since $z_{HV}(x^*)>0$ as shown above,
a.a.s.\ $\orH$ has a non-empty $(w,k+1)$-core $\coreH$. Recall that
$n$ and $m_{h-j}$ denote the number of vertices and hyperedges of
size $h-j$ in $\coreH$. Then a.a.s.\ the number of vertices in
$\coreH$ is $\bar n z_{HV}(x^*)+o(\bar n)$, and the number of
hyperedges of size $h-j$ in $\coreH$ is $\bar n
z_{H,h-j}(x^*)/(h-j)+o(\bar n)$. Since $z_{HV}(x^*)>0$ and
$z_{H,h-j}(x^*)>0$, we have a.a.s.\ $n\sim \alpha\bar n$ and
$m_{h-j}\sim\beta_{h-j}\bar n$, where $\alpha=z_{HV}(x^*)$ and
$\beta_{h-j}=z_{H,h-j}(x^*)/(h-j)$.

This proves the assertions about  $\orH\in \mathcal{M}_{\bar n,\bar
m,h}$. Lemma~\ref{uniformSimple} transfers them to
$\mathcal{G}_{\bar n,\bar m,h}$. \qed\ss

 There are several useful results that we will now derive recalling various pieces of the proof of Theorem~\ref{t:CoreDensity}. As noted at the start of that proof, the partition-allocation $g_{\tau}$ output by the RanCore
algorithm, if it is nonempty, is distributed as
$\mathcal{P}(V_{\tau},{\bf M}_{\tau},{\bf 0},k+1)$ conditional on
$V_{\tau}$ and ${\bf M}_{\tau}$. Let $n$ denote $|V_{\tau}|$ and
$m_{h-j}$ denote $|M_{\tau,h-j}|/(h-j)$ for all $0\le j\le w-1$.
Without loss of generality, by relabeling elements in $V_{\tau}$ and
${\bf M}_{\tau}$ in a canonical way, we can simplify the notation
$\mathcal{P}(V_{\tau},{\bf M}_{\tau},{\bf 0},k+1)$ to
$\mathcal{P}([n],{\bf M},{\bf 0},k+1)$, where ${\bf
M}=(M_{h-w+1},\ldots,M_h)$ and $M_i=[m_{i}]\times [i]$.
 The space $\mathcal{P}([n],{\bf
M},{\bf 0},k+1)$ is used in the proof of Theorem~\ref{t:ForCore} in
Section~\ref{probabilistic}.

\begin{lemma}\lab{l:exist}
 Assume $c_1k<h\bar m/\bar n<c_2k$ for some constants $c_2>c_1>1$.
Let $H$ be a random multihypergraph in $\mathcal{M}_{\bar n,\bar
m,h}$. Then, provided $k$ is sufficiently large, a.a.s.\ $H$ has a
nonempty $(w,k+1)$-core with average degree $O(k)$.
\end{lemma}
\proof Let $\bar\mu=h\bar m/\bar n$. Since $\bar\mu>c_1 k$ for some
$c_1>1$, the existence of a non-empty $(w,k+1)$-core  has
been shown in Theorem~\ref{t:de}.
 Let $x^*$ be as
defined in the statement of Theorem~\ref{t:CoreDensity} and let
$\delta=L_0/\bar n$. We have shown that $\delta=O(e^{-\Omega(k)})$
below~\eqn{eq:derivative} and $x^*\le 3\delta$ below~\eqn{eq:L}. Let
$z_B(x)$ and $z_{HV}(x)$ be defined the same as those functions
in~\eqn{eq:DL}--\eqn{eq:DefSystem} for $0\le x\le x^*$. Then clearly
$z_B(x^*)\le z_B(0)$ since $z_B'(x)\le -1$ for all $0\le x\le x^*$.
We also have $z_{HV}'(x)\ge -1/8$ for all $0\le x\le x^*$ when $k$
is large enough, as shown in the argument below~\eqn{eq:derivative}.
So $z_{HV}(x^*)\ge z_{HV}(0)-x^*/8$ for sufficiently large $k$.
Since $z_{HV}(0)=f_{k+1}(\bar \mu)=1-O(e^{-\Omega(k)})$ and
$x^*=O(e^{-\Omega(k)})$, we have $z_{HV}(x^*)=1-O(e^{-\Omega(k)})$.
Recall that $\mu(x)=(z_B(x)-z_L(x))/z_{HV}(x)$. Recall also that
$z'_B(x)-z'_L(x)\le 0$ by the argument below~\eqn{eq:derivative}.
Thus, we have $\mu(0)=O(k)$ since $h\bar m/\bar n<c_2 k$ and
$$
\mu(x^*)\le
\frac{z_B(0)-z_L(0)}{z_{HV}(x^*)}=\frac{z_B(0)-z_L(0)}{z_{HV}(0)}(1+O(e^{-\Omega(k)}))=O(k).
$$
By Theorem~\ref{t:CoreDensity}, the average degree of the
$(w,k+1)$-core of $H$ is asymptotically $\mu(x^*)$, which is bounded
by $O(k)$. \qed
\ss

 The following lemma gives a lower bound on the size of
the $(w,k+1)$-core of a random $h$-multihypergraph.
\begin{lemma}\lab{l:size}
 Assume $c_1k<h\bar m/\bar n<c_2k$ for some constants $c_2>c_1>1$.
Let $H$ be a random multihypergraph in $\mathcal{M}_{\bar n,\bar
m,h}$. Then a.a.s.\ the number of vertices in the $(w,k+1)$-core of
$H$ is $(1-O(e^{-\Omega(k)}))\bar n$.
\end{lemma}
\proof Let $n$ denote the number of vertices in the $(w,k+1)$-core
of $H$.  We showed just after~\eqn{eq:L}  that $x^*<3\delta$ where
we had $\delta=\exp(-\Omega_{c}(k))$. Since in each step at most $h$
heavy bins can disappear, the result follows
from~\eqn{tauvalue}.\qed \ss

We need the following lemma before proving
Theorem~\ref{t2:CoreDensity}.
\begin{lemma}\lab{l:denser}
 Assume $c_1k<h\bar m/\bar n<c_2k$ for some constants $c_2>c_1>1$. Let $\eps>0$ be fixed. Let $H_1$ be a random
 multihypergraph in
$\mathcal{M}_{\bar n,\bar m,h}$ and $H_2\in\mathcal{M}_{\bar n,\bar
m+\eps\bar n,h}$. Let $n_1$ and $n_2$ be the number of vertices in
the $(w,k+1)$-core of $H_1$ and $H_2$ respectively. Then a.a.s.\ we
have $|n_1-n_2|=O(e^{-\Omega(k)}\eps\bar n)$.
\end{lemma}
\proof Since $c_1k<h\bar m/\bar n<c_2 k$, by Lemma~\ref{l:exist},
the $(w,k+1)$-core $\widehat H_1$ of $H_1$ exists and the average
degree of $\widehat H_1$ is $O(k)$.
 Let $H_2$ be a random uniform multihypergraph obtained from
$H_1\cup \mathcal{E}$, where $\mathcal{E}$ is a set of $\eps\bar n$
hyperedges, each of which is a multiset of $h$ vertices, each of
which u.a.r.\ chosen from $[\bar n]$. Then $H_2\in\mathcal{M}_{\bar
n,\bar m+\eps\bar n,h}$. We say that the hyperedges in $\mathcal{E}$
are {\em marked}, and the other hyperedges in $H_2$ are {\em
unmarked}. Define a random process $(H_t^{(1)},H_t^{(2)})_{t \ge 0}$
as follows.

\begin{description}

\item{(i)} The process starts with $(H_0^{(1)},H_0^{(2)})=(H_1,H_2)$.

\item{(ii)} The RanCore algorithm is applied to $H_t^{(2)}$ for every
$t\ge 0$. The process $(H_t^{(1)},H_t^{(2)})_{t \ge 0}$ stops when
the RanCore algorithm running on $(H_t^{(2)})_{t\ge 0}$ terminates.

\item{(iii)} For every $t\ge 0$, if a marked hyperedge $x$ in
$H_{t-1}^{(2)}$ is updated to $x'$, then $x'$ remains marked in
$H_{t}^{(2)}$ and  $H_t^{(1)}$ is defined as $H_{t-1}^{(1)}$; if a
marked hyperedge $x$ is removed, also let $H_t^{(1)}=H_{t-1}^{(1)}$.

\item{(iv)} For every $t\ge 0$, if an unmarked hyperedge $x$ in
$H_{t-1}^{(2)}$ is updated or removed, do the same operation to $x$
in $H_{t-1}^{(1)}$ and define $H_{t}^{(1)}$ to be the resulting
hypergraph.
\end{description}

We call the random process $(H_t^{(i)})_{t\ge 0}$ for $i=1,2$
generated by $(H_t^{(1)},H_t^{(2)})_{t \ge 0}$ the {\em
$H_i$-process}. Note that the $H_1$-process is not equivalent to
running the RanCore algorithm on $H_1$, since the light balls are
not chosen u.a.r.\ in each step.

Instead of analysing $(H_t^{(1)},H_t^{(2)})_{t \ge 0}$ directly, we
consider $(g_t^{(1)},g_t^{(2)})_{t \ge 0}$, the corresponding
process obtained by considering the pairing-allocation model. Recall
that $H_1$ can be represented as dropping $h\bar m$ unmarked balls
u.a.r.\ into $\bar n$ bins with balls evenly partitioned into $\bar
m$ groups randomly and $H_2$ can be represented as dropping
$h\eps\bar n$ partitioned marked balls into $H_1$. The
partition-allocation $g_0^{(i)}$ for $i=1,2$ is obtained by putting
all balls contained in light bins of $H_i$ into one light bin.
Define $L_t^{(i)}$, $HV_t^{(i)}$, ${\bf m}_t^{(i)}$ and ${\bf
L}_t^{(i)}$, etc., for $i=1,2$ and for $t\ge 0$, the same way as in
the proof of Theorem~\ref{t:CoreDensity}, for the $H_i$-process.
Conditional on
 $L_0^{(i)}$, $V_{0}^{(i)}$, ${\bf M}_{0}^{(i)}$ and
${\bf L}_0^{(i)}$, $g_0^{(i)}$ is distributed as
$\mathcal{P}(V_{0}^{(i)},{\bf M}_{0}^{(i)},{\bf L}_0^{(i)},k+1)$ for
$i=1,2$ and all balls in $g_0^{(1)}$ are unmarked.

Let $\bar\mu$ denote the average degree of $H_1$ and let $\tau$ be
the time the $H_2$-process terminates. It is easy to show that
$g_{\tau}^{(1)}$ is distributed as $\mathcal{P}(V_{\tau}^{(1)},{\bf
M}_{\tau}^{(1)},{\bf L}_{\tau}^{(1)},k+1)$ conditional on the values
of $V_{\tau}^{(1)}$, ${\bf M}_{\tau}^{(1)}$ and ${\bf
L}_{\tau}^{(1)}$, since whenever a light ball is chosen, even not
uniformly at random, it results in recolouring or removal of heavy
balls that are uniformly chosen at random. We will later let the
RanCore algorithm be run
 on $g_{\tau}^{(1)}$ in the following steps and  apply the d.e.\ method to analyse the asymptotic behavior of
 this process.

First we show that $\tau=O(e^{-k}\bar n)$. The solution of the
differential equation system~\eqn{eq:DL}--\eqn{eq:DefSystem} tells
the asymptotic value of $L_t^{(2)}$ in every step $t$. Let
$x_{(2)}^{*}$ be the smallest root of $z_L^{(2)}(x)=0$. Since
$z_L^{(2)}(0)=O(e^{-\Omega(k)})$ and by the argument
below~\eqn{eq:derivative}, $z_L'^{(2)}(x)<-1/2$ for all $0\le
x<x_{(2)}^{*}$ provided $k$ sufficiently large,  we have
$x_{(2)}^{*}=O(e^{-\Omega(k)})$ and so $\tau=O(e^{-\Omega(k)}\bar
n)$.

Next we show that $n_2-n_1=O(e^{-\Omega(k)}\eps\bar n)$, assuming
the following three statements.
\begin{description}

\item{(S1)} The number of balls that are unmarked and light in $g_0^{(1)}$ but not
in $g_0^{(2)}$
 is bounded by $O(e^{-\Omega(k)}\eps\bar n)$.

\item{(S2)} The number of bins that begin heavy in the  $H_1$-process and become light in that process but remain heavy in the $H_2$-process up
to step $\tau$ is $O(e^{-\Omega(k)}\eps\bar n)$.

\item{(S3)} $L_{\tau}^{(1)}=O(e^{-\Omega(k)}\eps\bar n)$.

\end{description}

Run the Rancore algorithm on $g_{\tau}^{(1)}$. The differential
equation system~\eqn{eq:DL}--\eqn{eq:DefSystem} tells
 the asymptotic values of the
various random variables in $g_t^{(1)}$ for all $t\ge \tau$. Let
$x_{(1)}^{*}$ be the smallest positive root of $z_L^{(1)}(x)=0$.
Since $z_L^{(1)}(\tau/\bar n)=O(e^{-k}\eps)$ by (S3) and by the
argument below~\eqn{eq:derivative}, $z_L'^{(1)}(x)\le -1/2$ for all
$\tau/\bar n\le x<x_{(1)}^{*} $ provided $k$ sufficiently large, we
have $x_{(1)}^{*}-\tau/\bar n=O(e^{-k}\eps )$. We also have $-1/8\le
z_{HV}'(x)\le 0$ for sufficiently large $k$ for all $\tau/\bar n\le
x<x_{(1)}^{*} $ as explained in Lemma~\ref{l:exist}. So
$HV^{(1)}_{\tau}-n_1=O(e^{-k}\eps \bar n)$. Since
$n_2-HV^{(1)}_{\tau}$ counts the number of bins that are, or become
light in the $H_1$-process but stay heavy in the $H_2$-process, it
follows from (S1) and (S2)  that
$n_2-HV^{(1)}_{\tau}=O(e^{-k}\eps\bar n)$. So
$|n_1-n_2|=O(e^{-k}\eps\bar n)$.

It only remains to prove (S1)--(S3). We first show that (S3) follows
directly from (S1) and (S2).  $L_{\tau}^{(1)}$ counts two types of
light balls. The first type comes from balls that are unmarked and
light in $g_0^{(1)}$ but not in $g_0^{(2)}$. By (S1), the number of
these balls is a.a.s.\ $O(e^{-\Omega(k)}\eps\bar n)$. The second
type comes from balls that begin heavy and become light  in the
$H_1$-process but stay heavy in the $H_2$-process. By (S2), the
number of these balls is a.a.s.\ $k\cdot O(e^{-\Omega(k)}\eps\bar
n)=O(e^{-\Omega(k)}\eps\bar n)$. Thereby (S3) follows.

Next we show (S1).
 At step 0, clearly the set of unmarked light balls in $g_0^{(2)}$
is a subset of those in $g_0^{(1)}$.
 The number of light balls in $g_0^{(1)}$ is a.a.s.\ $(1-f_k(\bar\mu))\bar\mu\bar
 n=O(e^{-\Omega(k)}\bar n)$ as shown in the proof of
 Theorem~\ref{t:CoreDensity} and hence the number of light vertices of $H_1$
 is a.a.s.\
  $O(e^{-\Omega(k)}\bar n)$. Since each multihyperedges in $\mathcal{E}$ is a random
  multihyperedges, the expected number of those which contains a
  light vertex in $H_1$ is $O(e^{-\Omega(k)}\eps\bar n)$, hence the
  number of light vertex in $H_1$ that become heavy after the
  hyperedges in $\mathcal{E}$ being dropped is a.a.s.\ $O(e^{-\Omega(k)}\eps\bar
  n)$ and each of these vertex/bin contains at most $k$ unmarked balls. Thus (S1) follows.

Now we show (S2). Recall that $H_1$ is represented as dropping
$h\bar m$ unmarked balls u.a.r.\ into $\bar n$ bins and $H_2$ is
obtained by dropping $h\eps \bar n$ extra marked balls u.a.r.\ into
the $\bar n$ bins in $H_1$. Recall that $\coreH_1$ denotes the
$(w,k+1)$-core of $H_1$. The number of bins that begin heavy in the
$H_1$-process and become light in that process but remain heavy in
the $H_2$-process up to step $\tau$ is at most the number of
bins/vertices not in $\coreH_1$ which receive at least one marked
balls after dropping $h\eps\bar n$ marked balls u.a.r.\ into the
$\bar n$ bins. By Lemma~\ref{l:size}, the number of vertices/bins in
$\coreH_1$ is a.a.s.\ $(1-O(e^{-\Omega(k)}))\bar n$. Then for each
marked ball, the probability that it is dropped into a bin not in
$\coreH_1$ is $O(e^{-\Omega(k)})$. By Lemma~\ref{l:super}, the
number of marked balls dropped into bins not in $\coreH_1$ is
a.a.s.\ $O(e^{-\Omega(k)}\eps\bar n)$. Hence the number of bins that
are not in $\coreH_1$ and receive at least one marked balls is
a.a.s.\ $O(e^{-\Omega(k)}\eps\bar n)$.\qed\ss

\no {\bf Proof of Theorem~\ref{t2:CoreDensity}.\ } Let $H_1$ be a
random uniform multihypergraph with  average degree $\bar\mu$ and
let $H_2$ be a random uniform multihypergraph obtained from $H_1\cup
\mathcal{E}$, where $\mathcal{E}$ is a set of $\eps\bar n$
hyperedges, each of which is a multiset of $h$ vertices, each of
which is uniformly chosen from $[\bar n]$.

For $i=1,2$, let $\widehat H_i$ be the $(w,k+1)$-core of $H_i$ and
let $m_{h-j}^{(i)}$ be the number of hyperedges with size $h-j$ in
$\widehat H_i$. We first show that
$$
\sum_{j=0}^{w-1} (w-j)m_{h-j}^{(2)}-\sum_{j=0}^{w-1}
(w-j)m_{h-j}^{(1)}\ge w\eps \bar n/2.
$$

Clearly $\widehat H_1$ is a subgraph of $\widehat H_2$. Let $n_i$
denote the number of vertices in $\widehat H_i$ and let $[n_i]$
denote the set of vertices in $\widehat H_i$. By Lemma~\ref{l:size},
a.a.s.\ $n_1=(1-O(e^{-\Omega(k)}))\bar n$. Then for any hyperedge
$x\in\mathcal{E}$, the probability that all vertices in $x$ are
contained in $[n_1]$ is $1-O(e^{-\Omega(k)})$. So the expected
number of hyperedges in $\mathcal{E}$ lying completely in $[n_1]$ is
$(1-O(e^{-\Omega(k)}))\eps\bar n$. By the Chernoff bound, originally
given in~\cite[Theorem 1]{C3}, we have a.a.s.\ the number of
hyperedges in $\mathcal{E}$ lying completely in $[n_1]$ is at least
$\eps\bar n/2$ for sufficiently large $k$. So it follows immediately
that a.a.s.,
$$
\sum_{j=0}^{w-1} (w-j)m_{h-j}^{(2)}-\sum_{j=0}^{w-1}
(w-j)m_{h-j}^{(1)}\ge w\eps \bar n/2.
$$
For simplicity, let $S(i)$ denote $\sum_{j=0}^{w-1}
(w-j)m_{h-j}^{(i)}$ for $i=1,2$. Recall that $\kappa(\widehat H_i)$
denotes $S(i)/n_i$. Then a.a.s.,
\begin{eqnarray*}
\kappa(\widehat H_2)-\kappa(\widehat
H_1)=\frac{S(2)}{n_2}-\frac{S(1)}{n_1}\ge \frac{(S(1)+w\eps\bar
n/2)-S(1)\cdot n_2/n_1}{n_2}.
\end{eqnarray*}
By Lemma~\ref{l:denser}, a.a.s.\ $n_2-n_1=O(e^{-\Omega(k)})\eps\bar
n$, i.e.\ $n_2/n_1-1\le f(k)\eps$ for some function
$f(k)=O(e^{-\Omega(k)})$. Then a.a.s.,
\begin{eqnarray}
\kappa(\widehat H_2)-\kappa(\widehat H_1)\ge \frac{w\eps\bar
n/2-O(f(k)\eps S(1))}{n_2}\ge w\eps/4>0,\lab{eq:monotone}
\end{eqnarray}
for sufficiently large $k$ and for every $\eps>0$, since
$S(1)=O(k)\bar n$ and $n_2=(1-O(e^{-\Omega(k)}))\bar n$.

 By
Theorem~\ref{t:CoreDensity}, for given $h>w>0$ and sufficiently
large $k$, a.a.s.\ $\kappa(\coreH)= c(\bar\mu)+o(1)$, where
$c(\bar\mu)$ is a constant depending only on $\bar\mu$.
 The inequality~\eqn{eq:monotone} implies that $c(\bar\mu)$ is an
increasing function of $\bar\mu$.   \qed\ss

\no {\bf Proof of Corollary~\ref{c:CoreDensity}.\ } By
Theorem~\ref{t2:CoreDensity}, there exists a unique critical value
of $\bar\mu$ such that a.a.s.\ $\kappa(\coreH)=k+o(1)$ and so there
exists a threshold function $\bar m=f(\bar n)$ of $\mathcal{M}_{\bar
n,\bar m,h}$ for the graph property $\mathcal{T}$. Then this holds
as well in $\mathcal{G}_{\bar n,\bar m,h}$ by
Lemma~\ref{uniformSimple}.\qed\ss

The differential equations in Theorem~\ref{t:CoreDensity} are only
used in a theoretical way to show properties of the $(w,k+1)$-core,
and we do not have an analytic solution. However, they can be
numerically solved when the values of $h$, $w$, $k$ and $\mu$ are
given.  Table~\ref{t:computation} gives the results of some
computations, where $h$, $w$ and $k$ are given,  $\widetilde\mu$
denotes the expected average degree of the hypergraph $\orH$ at the
threshold for $\cal T$ given in Corollary~\ref{c:CoreDensity}, and
$\widehat\mu$ denotes the corresponding average degree of its core
$\coreH$. Even though our results on the concentration of the size
and density of the $(w,k+1)$-core and the threshold of property
$\cal T$ only cover for the case of sufficiently large $k$, our
numerical computation results as shown in the table do coincide with
our simulation results. Hence we believe that
Theorem~\ref{t:CoreDensity},~\ref{t2:CoreDensity} and
Corollary~\ref{c:CoreDensity} actually hold for all $k\ge 1$. By
Corollary~\ref{c2:lbalancing}, discussed in the next section,
$\widetilde\mu$ is also our main target, the threshold for
orientability. Note that $\widehat\mu$ must be at least $hk/w$ by
the definition of property $\cal T$, and that it follows from the
trivial upper bound of the orientability threshold given in the
introduction part that $\widetilde \mu$ is at most $hk/w$.
\begin{table}[htb]
\begin{center}
 \vspace{.2cm}

\begin{tabular}[b]{|c|c|c|c|c|}    \hline
    $h$&$w$&$k$&$\widetilde\mu$&$\widehat\mu$\\ \hline
    $3$&$2$&$4$&$5.485$&$6.65086$\\ \hline
    $3$&$2$&$10$&$14.766$&$15.5872$\\ \hline
    $3$&$2$&$40$&$59.991$&$60.0773$\\ \hline
    $10$&$2$&$4$&$19.99999$&$20.0003$\\ \hline
   \end{tabular}
\vspace{.3cm} \lab{t:computation}
 \caption{Some numerical computation results}
\end{center}
\end{table}

\section{The $(w,k)$-orientability of the $(w,k+1)$-core }
\lab{orientability}

In this section we prove Corollary~\ref{c2:lbalancing} assuming
Theorem~\ref{t:ForCore}, and study the basic network flow
formulation of the problem that is used in the next section to prove
Theorem~\ref{t:ForCore}.

The following lemma is in preparation for proving that $\coreH$, the
$(w,k+1)$-core of $\orH\in\mathcal{M}_{\bar n,\bar m,h}$, if not
empty, a.a.s.\ has property $\mathcal{A}(\gamma)$ for some
$0<\gamma<1$.
\begin{lemma}\lab{l:SmallSet}
Let $\orH\in\mathcal{M}_{\bar n,\bar m,h}$ and let $\coreH$ be the
$(w,k+1)$-core of $\orH$. Let $c_1>1$  be a constant that can depend
on $k$.  Then there exists a constant $0<\gamma=\varphi(k,c_1)$
depending only on $k$ and $c_1$, such that a.a.s.\ there exists no
$S\subset V(H)$ with $|S|<\gamma \bar n$ and at least $c_1|S|$
hyperedges partially contained in $S$. More specifically, when
$c_1\ge 2$ and $c_1<h^2e^2\bar\mu$, we may choose
$\gamma=\varphi(k,c_1)=\left(c_1/h^2e^2\bar\mu\right)^{2}$.
\end{lemma}
\proof Let $s$ be any integer such that $0<s<n$ and let $r=s/n$. Let
$Y$ denote the number of $S$ with $|S|= s$ and at least $c_1s$
hyperedges partially contained in $S$.  The probability for a given
hyperedge to be partially contained in $S$ is at most $ \binom{h}{2}
(s/\bar n)^2<h^2 r^2$. Then the probability that there are at least
$c_1s$ such hyperedges is at most
$$
\binom{\bar m}{c_1s}(h r)^{2c_1s}.
$$
Since there are $\binom{\bar n}{s}$ ways to choose $S$,
\begin{eqnarray*}
\ex(Y)&=&\displaystyle\sum_{s\le \gamma\bar n}\binom{\bar n}{s}\binom{\bar m}{c_1 s}(hr)^{2c_1s}\\
&\le& \sum_{\ln \bar n\le s\le \gamma\bar n}\left(\frac{e\bar
n}{s}\right)^s\left(\frac{e\bar m}{c_1s}\right)^{c_1s}(hr)^{2c_1s}+\sum_{1\le s\le \ln \bar n}\bar n^s\bar m^{c_1s}
\left(\frac{hs}{\bar n}\right)^{2c_1s}\\
&=& \sum_{\ln \bar n\le s\le \gamma\bar n}\left(h^{2c_1}e^{1+c_1}r^{c_1-1}\left(\frac{\bar\mu}{c_1}\right)^{c_1}\right)^{s}
+\sum_{1\le s\le \ln\bar n}\left(\frac{(\bar\mu h^2s^2)^{c_1}}{\bar n^{c_1-1}}\right)^s\\
&\le& \sum_{\ln\bar n\le s\le \gamma \bar n}\left(\bar C
r^{c_1-1}\bar\mu^{c_1}\right)^{s} +\ln\bar n\cdot\frac{(\bar\mu h^2 \ln^2 \bar n)^{c_1}}{\bar n^{c_1-1}}\\
&=& \sum_{\ln\bar n\le s\le \gamma\bar n}\left(\bar C
r^{c_1-1}\bar\mu^{c_1}\right)^{s}+o(1),
\end{eqnarray*}
for some constant $0<\bar C=\bar C(c_1)\le
\left(h^2/c_1\right)^{c_1}e^{c_1+1}$. Choose
$$
\gamma<\left(\frac{c_1}{h^2e\bar\mu}\right)^{\frac{c_1}{c_1-1}}e^{-\frac{1}{c_1-1}}.
$$
Then $\bar C \gamma^{c_1-1}\bar\mu^{c_1}<1$. So there exist
$0<\beta<1$, such that $\bar C \gamma^{c_1-1}\bar\mu^{c_1}<\beta$,
for all $r\le \gamma$. When $c_1\ge 2$ and  $c_1/h^2e^2\bar\mu<1$,
$$
\left(\frac{c_1}{h^2e\bar\mu}\right)^{\frac{c_1}{c_1-1}}e^{-\frac{1}{c_1-1}}>\left(\frac{c_1}{h^2e^2\bar\mu}\right)^{c_1/(c_1-1)}
>\left(\frac{c_1}{h^2e^2\bar\mu}\right)^2.
$$
 Hence we may simply choose
$\gamma=(c_1/h^2e^2\bar\mu)^{2}$. Then
$$
\sum_{\ln\bar n\le s\le \gamma\bar n}\left(\bar C
r^{c_1-1}\bar\mu^{c_1}\right)^{s}<\sum_{\ln\bar n\le s\le \gamma\bar
n}\beta^s=O(\beta^{\ln\bar n})=o(1).
$$
Hence we have $\ex(Y)=o(1)$.\qed\ss

The following corollary shows that the same property is shared by
$\coreH$.
\begin{cor}\lab{c:SmallSet}
 Let $\orH\in\mathcal{M}_{\bar n,\bar m,h}$ and let $\coreH$ be the $(w,k+1)$-core of $\orH$.
Let  $c_1$ be a constant that can depend on $k$, with the constraint
that $2\le c_1<h^2e^2\bar\mu$.
Let $0<\gamma=\varphi(k,c_1)=\left(c_1/h^2e^2\bar\mu\right)^{2}$.
Then a.a.s.\ for all $S\subset V(\coreH)$ with $|S|<\gamma n$, the
number of hyperedges partially contained in $S$ is less than
$c_1|S|$.
\end{cor}
\proof
Let $n$ be the number of vertices in $\coreH$ and $D$ the sum of
degrees of vertices in $\coreH$. For any hyperedge $x\in \coreH$,
let $x^{+}$ denote its corresponding hyperedge in $\orH$. Obviously
$n\le\bar n$. Combining with Lemma~\ref{l:SmallSet} and the fact
that for any $S\subset V(\coreH)$, a hyperedge $x$ is partially
contained in $S$ only if $x^{+}$ is partially contained $S$ in
$\orH$, Corollary~\ref{c:SmallSet} follows.\qed\ss

We next show that $\coreH$, if not empty, a.a.s.\ has property
$\mathcal{A}(\gamma)$, defined in Definition~\ref{d:A}, for some
certain value of $\gamma$.
\begin{cor}\lab{c:SmallSet2} Assume that $\orH\in\mathcal{M}_{\bar n,\bar m,h}$ with $\bar m\le hk/w$, and that $\coreH$ is the $(w,k+1)$-core of $\orH$. Let $\gamma=e^{-4}h^{-6}/4$. Then provided $k\ge 4w$, a.a.s.\ either $\coreH$ is empty or
$\coreH$ has property $\mathcal{A}(\gamma)$.
\end{cor}
\proof Apply Lemma~\ref{c:SmallSet} with $c_1=k/2w$. Clearly
$c_1<ch^2e^2k$, and $c_1\ge 2$ provided $k\ge 4w$.  Then $\gamma\le
\phi(k,c_1)$. By Definition~\ref{d:A}, $\coreH$ a.a.s.\ has property
$\mathcal{A}(\gamma)$. \qed \ss

\no {\bf Proof of Corollary~\ref{c2:lbalancing}} Let $\coreH$ be the
$(k+1)$-core of the random multihypergraph $H\in \mathcal{M}_{\bar
n, \bar m,h}$. Let $\eps>0$ be any constant. By
Theorem~\ref{t2:CoreDensity}, there exists a constant $\delta>0$,
such that a.a.s.\ if $\bar m\le f(\bar m)-\eps\bar n$, then
$\sum_{j=0}^{w-1}(w-j)m_{h-j}\le kn-\delta n$. By
Theorem~\ref{t:ForCore} and Corollary~\ref{c:SmallSet2},
  there exists a constant $N$ depending only
on $h$ and $w$ such that provided $k>N$, $\coreH$ a.a.s.\ has a
$(w,k)$-orientation. On the other hand, if $\bar m\ge f(\bar
m)+\eps\bar n$, then a.a.s.\ $\sum_{j=0}^{w-1}(w-j)m_{h-j}\ge
kn+\delta n$, and hence clearly $\coreH$ is not $(w,k)$-orientable.
Therefore $f(\bar n)$ is a sharp threshold function for the
$(w,k)$-orientation of $\mathcal{M}_{\bar n, \bar m,h}$. By
Lemma~\ref{uniformSimple}, $f(\bar n)$ is also a sharp threshold
function for the $(w,k)$-orientation of $\mathcal{G}_{\bar n, \bar
m,h}$. \qed \ss

Let $\core2H$ be a non-uniform multihypergraph with the sizes of
hyperedges between $h-w+1$ to $h$. In the rest of the chapter, we
will use the following notations. Let $E_{h-j}:=\{x\in E(\core2H):
|x|=h-j\}$.  For any given $S\subset [n]$, let $
m_{h-j,i}(S):=|\{x\in E_{h-j}: |x\cap S|=i\}|$ for any $0\le i\le
h-j$.  When the context is clear of which set $S$ is referred to, we
may drop $S$ from the notation. Let $\barS$ denote the set
$[n]\setminus S$ and let $d(S)$ denote the sum of degrees of
vertices in $S$.

Recall from above the statement of
Corollary~\ref{c:lbalancing} in Section~\ref{mainResult} that for
any $S\subset V(\core2H)$, $\core2H_S$ denotes the subgraph
$w$-induced by $S$. The following Lemma generalises Hakimi's
theorem~\cite[Theorem 4]{H2} for graphs. It is proved using network
flow and the max-flow min-cut theorem, along the lines of the
standard  techniques discussed in~\cite{CCPS,S2}. This setting was
used before in connection with the load balancing problem
in~\cite[Section 3.3]{SEK}.
\begin{lemma}\lab{l:subgraph}
A multihypergraph $\core2H$ with sizes of hyperedges between $h-w+1$
and $h$ has a $(w,k)$-orientation if and only if
\remove{d(\orH_S)-(h-w)e(\orH_S)\le k|S|,} $\kappa(\core2H_S)\le k$
for all $S\subset V(\core2H)$.
\end{lemma}
\proof  Formulate a network flow problem on a network $G^*$ as
follows. Let $L$ be a set of vertices, each of which represents a
hyperedge of $\core2H$, and $R$ be a set of $n$ vertices, each of
which represents a vertex in $\core2H$. For any $u\in L$, and $v\in
R$, $uv$ is an edge in $G^*$ if and only if $v\in u$ in $\core2H$.
Add vertices $a$ and $b$ to $G^*$, such that $a$ is linked to every
vertex in $L$, and $b$ is linked to every vertex in $R$.  Let $c:E(
G^*)\rightarrow {\bf N}^+$ be defined as $c(au)=w-j$ for every $u\in
L$ such that the degree of $u$ is $h-j$, $c(vb)=k$ for every $v\in
R$, and $c(uv)=1$ for every $uv\in E(G^*)$. Then $\core2H$ has a
$(w,k)$-orientation if and only if $ G^*$ has a flow of size
$\sum_{j=0}^{w-1}(w-j)m_{h-j}$ from $a$ to $b$. By the max-flow
min-cut Theorem, $ G^*$ has a flow with all edges incident with $a$
saturated if and only if
\begin{equation}\lab{eq:Cut2}
c(\delta(C))\ge\sum_{j=0}^{w-1}(w-j)m_{h-j}, \ \  \mbox{ for all
(a,b)-cuts}\ C.
\end{equation}

  \begin{figure}[htb]
\vbox{\vskip .8cm
 \hbox{\centerline{\includegraphics[width=8cm]{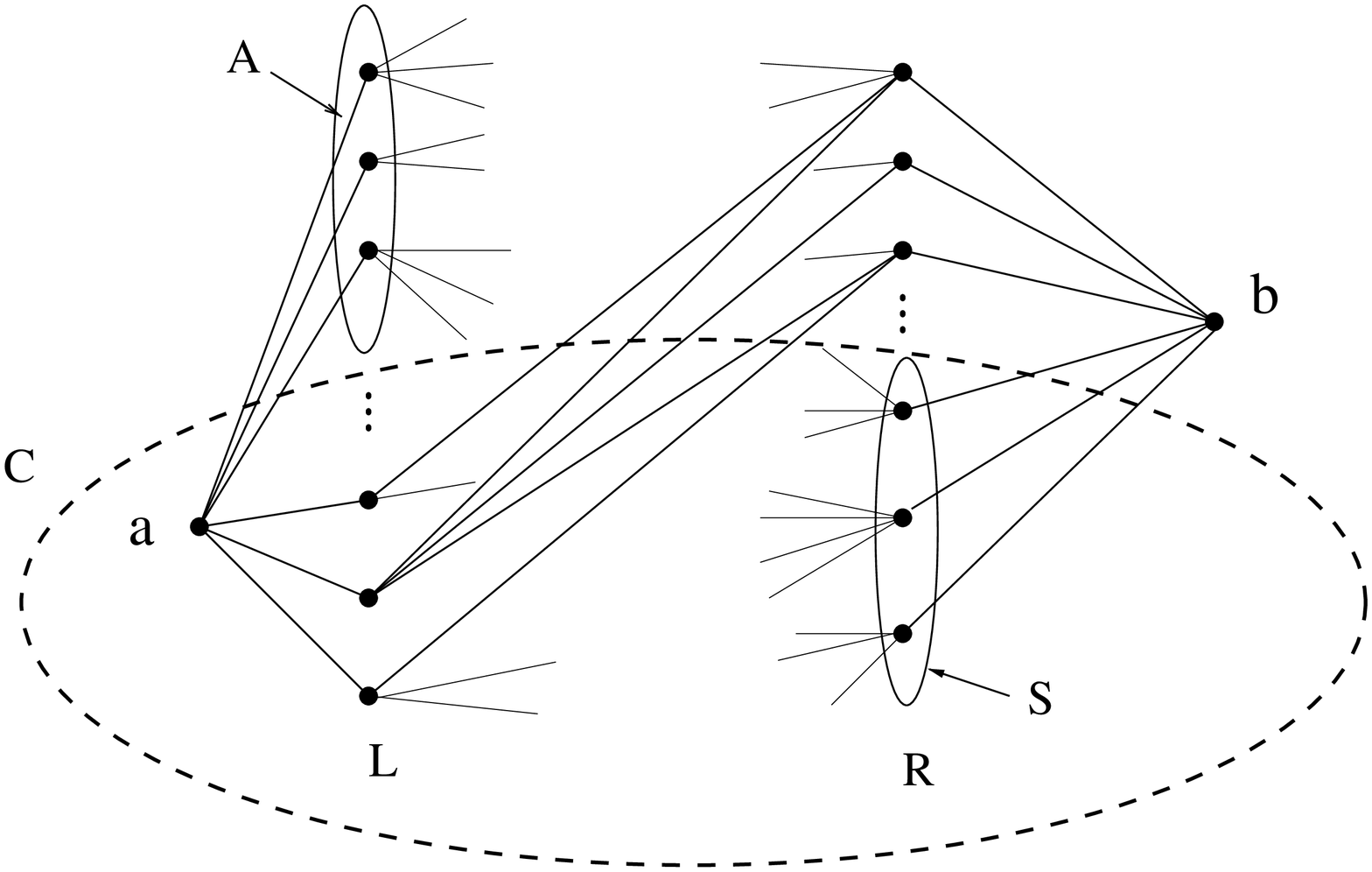}}}

\vskip .5cm \smallskip} \caption{\it  A cut $C$ in the graph $G^*$}

\lab{f:flow}

\end{figure}

As an example in Figure~\ref{f:flow}, $A\subset L$ is a set of
hyperedges in $\core2H$, and $S\subset R$ is a set of vertices in
$\core2H$. Let $C=\{a\}\cup \barA\cup S $ define a cut of $G^*$.
Then the condition in~\eqn{eq:Cut2} is equivalent to
\begin{equation}\lab{eq:Cut}
\forall C ,\ \ c(\delta(C))=k|S|+\sum_{j=0}^{w-1}\left(\sum_{x\in
A\cap E_{h-j}}(w-j)+\sum_{x\in E_{h-j}\setminus
A}|x\cap\barS|\right)\ge \sum_{j=0}^{w-1}(w-j)m_{h-j},
\end{equation}
Let $A^*:=\{x\in E_{h-j}: |x\cap S|\le h-w\}$. Clearly $A^*$
minimizes $c(\delta(C))$ for a given $S$. Therefore we only need to
check~\eqn{eq:Cut} when $A=A^*$. The condition in~\eqn{eq:Cut} is
then equivalent to
\begin{eqnarray*}
\sum_{j=0}^{w-1}\left(\sum_{x\in  E_{h-j}\setminus
A}(w-j)-\sum_{x\in E_{h-j}\setminus A}|x\cap \barS|\right)\le k|S|.
\end{eqnarray*}
For any hypergraph $\core2H$, let $\beta(\core2H)$ denote
the number of hyperedges in $\core2H$. Recall from the statement
above Theorem~\ref{t2:CoreDensity} that
$|S|\kappa(\core2H_S)=d(\core2H)-(h-w)\beta(\core2H)$. Since
\begin{eqnarray}
&&\sum_{j=0}^{w-1}\left(\sum_{x\in  E_{h-j}\setminus
A}(w-j)-\sum_{x\in E_{h-j}\setminus A}|x\cap
\barS|\right)=\sum_{j=0}^{w-1}\sum_{x\in  E_{h-j}\setminus
A}(w-j)-(h-j-|x\cap S|)\nonumber\\
&&\hspace{0.8cm}=\sum_{x\notin A}|x\cap S|-\sum_{x\notin
A}(h-w)=d(\core2H_S)-(h-w)\beta(\core2H_S)=|S|\kappa(\core2H_S),\lab{eq:sub}
\end{eqnarray}
Lemma~\ref{l:subgraph} follows. \qed \ss

The next corollary follows immediately.
\begin{cor}\lab{c:subgraph}
A hypergraph $H$ in $\mathcal{G}_{\bar n,\bar m,h}$ has a
$(w,k)$-orientation if and only if for every $S\subset V(H)$,
 $ \remove{d(\orH_S)-(h-w)e(\orH_S)\le k|S|.}
\kappa(\orH_S)\le k.$
\end{cor}
\no{\bf Proof of Corollary~\ref{c:lbalancing}.\ } This follows
directly from Corollary~\ref{c2:lbalancing} and
Corollary~\ref{c:subgraph}.\qed \ss

For any vertex set $S$, define
\begin{equation}
\partial^*(S)=d(S)-\sum_{j=0}^{w-1}\sum_{i=w-j+1}^{h-j}(i-(w-j))m_{h-j,i},\lab{eq:expansion}
\end{equation}
which measures a type of  expansion in the hypergraph.  For each
hyperedge $x$ of size $h-j$ which intersects $S$ with $i$ vertices,
 its contribution to $\partial^*(S)$ is $w-j\ge 0$ if $i\ge w-j+1$ and $i\ge 0$ otherwise. Therefore $\partial^*(S)\ge 0$ for any $S$.
\remove{ Recall that in a graph, the expansion rate of a vertex set
is defined as $\partial (S)/|S|$, where $\partial(S)$ is the number
of edges with exactly one end in $S$. It can be easily seen that
$\partial^*(S)$ is closely related to $\partial(S)$. For instance,
given the sum of degrees of vertices in $S$, both $\partial^*(S)$
and $\partial(S)$ are maximized when $S$ is an independent set and
minimized when $S$ is a component of the graph.} The following lemma
characterises the existence of the $(w,k)$-orientation of $\core2H$
in terms of $\partial^*(S)$.
\begin{lemma}\lab{l:expansion}
Let $\core2H$ be a multihypergraph whose hyperedges all have sizes
between $h-w+1$ and $h$ inclusively. Then the following two
properties of $\core2H$ are equivalent:
\begin{description}
\item{(i)} \remove{d(\orH_S)-(h-w)e(\orH_S)\le k|S|,} $\kappa(\core2H_S)\le k$, \ \ \ for all
$S\subset V(\core2H)$;
\item{(ii)} $\partial^*(S)\ge
k|S|+\left(\sum_{j=0}^{w-1}(w-j)m_{h-j}\right)-kn$, \ \ \ for all
$S\subset V(\core2H)$.
\end{description}
\end{lemma}
\proof  Let $\beta(\core2H_S)$ denote the number of
hyperedges in $\core2H_S$. We show that for any $S\subset
V(\core2H)$,
$\kappa(\core2H_S)|S|=d(\core2H_S)-(h-w)\beta(\core2H_S)\le k|S|$ if
and only if $\partial^*(\barS)\ge
k|\barS|+\Bigg(\sum_{j=0}^{w-1}(w-j)m_{h-j}\Bigg)-kn$. Then
Lemma~\ref{l:expansion} follows immediately. Note from the
definition of $A^*$, we have for any $x\in E_{h-j}\setminus A^*$,
$|x\cap S|\ge h-w+1$ and hence $|x\cap \barS|\le
(h-j)-(h-w+1)=w-j-1$. By~\eqn{eq:sub}, for any $S\subset
V(\core2H)$,
\begin{eqnarray*}
&&d(\core2H_S)-(h-w)\beta(\core2H_S)\le k|S|\\
&&\hspace{0.7cm}\Longleftrightarrow \sum_{j=0}^{w-1}\left(\sum_{x\in
E_{h-j}\setminus A^*}(w-j)-\sum_{x\in E_{h-j}\setminus A^*}|x\cap
\barS|\right)\le
kn-k|\barS|\\
&&\hspace{0.7cm}\Longleftrightarrow
\sum_{j=0}^{w-1}\sum_{i=0}^{w-j-1}
(w-j-i)m_{h-j,i}(\barS)\le kn-k|\barS|\\
&&\hspace{0.7cm}\Longleftrightarrow
\sum_{j=0}^{w-1}(w-j)m_{h-j}-\sum_{j=0}^{w-1}\left(\sum_{i=w-j}^{h-j}(w-j)m_{h-j,i}(\barS)+\sum_{i=0}^{w-j-1}im_{h-j,i}(\barS)\right)\le
kn-k|\barS|\\
&&\hspace{0.7cm}\Longleftrightarrow \partial^*(\barS)\ge
k|\barS|+\Bigg(\sum_{j=0}^{w-1}(w-j)m_{h-j}\Bigg)-kn.\qed
\end{eqnarray*}

It follows from Lemma~\ref{l:subgraph} and Lemma~\ref{l:expansion}
that $\core2H$ is $(w,k)$-orientable if and only if
Lemma~\ref{l:expansion} (ii) holds.

Without loss of generality, we assume
$\sum_{j=0}^{w-1}(w-j)m_{h-j}-kn\le 0$. Otherwise,
condition~\eqn{eq:Cut2} is violated by taking $C=\{a\}\cup L\cup R$.
The following lemma shows that, instead of checking conditions in
Lemma~\ref{l:expansion} (ii), we can check that certain other events
do not occur.

For any $S\subset V(\core2H)$, let
\begin{equation}\lab{eq:eta}
q_{h-j}(S)=\sum_{i=1}^{h-j} i m_{h-j,i},\ \ \
\eta(S)=\sum_{j=0}^{w-1}\sum_{i=1}^{h-j-1}m_{h-j,i}.
\end{equation}
In other words, $q_{h-j}(S)$ denotes the contribution to $d(S)$ from
hyperedges of size $h-j$ and $\eta(S)$ denotes the number of
hyperedges which intersect both $S$ and $\barS$. When the context is
clear, we may use $q_{h-j}$ and $\eta$ instead to simplify the
notation.

Recall that given a vertex set $S$, a hyperedge $x$ is partially
contained in $S$ if $|x\cap S|\ge 2$. Let $\rho(S)$ denote the
number of hyperedges partially contained in $S$ and let $\nu(S)$
denote the number of hyperedges intersecting $S$.
\begin{lemma}\lab{l:deterministic}

Suppose that for some $S\subset V(\core2H)$,
\begin{equation}
\partial^*(S)<
k|S|+\left(\sum_{j=0}^{w-1}(w-j)m_{h-j}\right)-kn.
\lab{eq:FlowCondition1}
\end{equation}
Then all of the following hold:
\begin{description}
\item{(i)} $\rho(\barS)>k|\barS|/w$;
\item{(ii)} $ \nu(S)<k|S|$;
\item{(iii)} $ (h-w)\rho(S)>d(S)-k|S|$;
\item{(iv) } 
 if, in addition,
$\displaystyle \sum_{j=0}^{w-1}\frac{w-j}{h-j}q_{h-j}(S)\ge (1-\delta)k|S|$
for some $\delta>0$, then $\eta(S)<h^2\delta k|S|$.
\end{description}
\end{lemma}


\proof Let $s$ and $\bar s$ denote $|S|$ and $|\barS|$ respectively.
If~\eqn{eq:FlowCondition1} is satisfied, then

$$
d(S)-\sum_{j=0}^{w-1}\sum_{i=w-j+1}^{h-j}(i-(w-j))m_{h-j,i}<ks-\left(kn-\sum_{j=0}^{w-1}(w-j)m_{h-j}\right)=\sum_{j=0}^{w-1}(w-j)m_{h-j}-k\bar
s.
$$
Hence
\begin{eqnarray}
k\bar s&<&\sum_{j=0}^{w-1}(w-j)m_{h-j}-\sum_{j=0}^{w-1}\sum_{i=1}^{h-j}im_{h-j,i}+\sum_{j=0}^{w-1}\sum_{i=w-j+1}^{h-j}(i-(w-j))m_{h-j,i}\nonumber\\
&=&\sum_{j=0}^{w-1}(w-j)m_{h-j,0}+\sum_{j=0}^{w-1}\sum_{i=1}^{w-1-j}(w-j-i)m_{h-j,i}\nonumber\\
&\le&
w\sum_{j=0}^{w-1}\left(m_{h-j,0}+\sum_{i=1}^{w-1-j}m_{h-j,i}\right).\nonumber
\end{eqnarray}
Since
$$
m_{h-j,0}= |\{x\in E_{h-j}: |x\cap \barS|=h-j\}|,
$$
and
$$
\sum_{i=1}^{w-1-j}m_{h-j,i}
\le w|\{x\in E_{h-j}: 2\le|x\cap \barS|\le h-j-1\}|,
$$
(this is because $1\le i\le w-1-j$ and so $h-j-i\le h-j-1$ and
$h-j-i\ge h-(w-1)\ge 2$), we have
$$
k\bar s<w|\{x\in E(G):|x\cap \barS|\ge 2\}|.
$$
This proves part (i). Again, if~\eqn{eq:FlowCondition1} is
satisfied, then
$$
\sum_{j=0}^{w-1}\sum_{i=w-j+1}^{h-j}(i-(w-j))m_{h-j,i}>d(S)-ks+\left(kn-\sum_{j=0}^{w-1}(w-j)m_{h-j}\right).
$$
Since
$$
\sum_{j=0}^{w-1}\sum_{i=w-j+1}^{h-j}(i-(w-j))m_{h-j,i}\le
\sum_{i=2}^h (i-1)|\{x:|x\cap S|=i\}|=d(S)-\nu(S),
$$
 we have
\begin{equation}
d(S)-\nu(S)>d(S)-ks+\left(kn-\sum_{j=0}^{w-1}(w-j)m_{h-j}\right).\nonumber
\end{equation}
Since $kn-\sum_{j=0}^{w-1}(w-j)m_{h-j}>0$, this directly leads to
part (ii). Since
$$
\sum_{j=0}^{w-1}\sum_{i=w-j+1}^{h-j}(i-(w-j))m_{h-j,i}\le
 (h-w)|\{x:|x
 \cap S|\ge 2\}|,
$$
we have
\begin{equation}
|\{x:|x\cap S|\ge
2\}|>d(S)-ks+\left(kn-\sum_{j=0}^{w-1}(w-j)m_{h-j}\right).\nonumber
\end{equation}
Since $kn-\sum_{j=0}^{w-1}(w-j)m_{h-j}>0$, this  proves part (iii).
Now we prove part (iv).
 Let
$t_{h-j}=1-(h-j)m_{h-j,h-j}/q_{h-j}$.
Note that $d(S)=\sum_{j=0}^{w-1}q_{h-j}$ and
$q_{h-j}=\sum_{i=1}^{h-j}im_{h-j,i}$. For each hyperedge $x$ of size
$h-j$ which intersects $S$ with $i$ vertices, its contribution to
$q_{h-j}$ (and thus to $\partial^*(S)$) is
\begin{itemize}
\item $i\cdot (w-j)/(h-j)$, if $i=h-j$;
\item $i\cdot (w-j)/i\ge i\cdot (w-j)/(h-j-1)$, if $w-j+1\le i\le
h-j-1$;
\item $i\ge i\cdot (w-j)/(h-j-1)$, if $1\le i\le w-j$;
\end{itemize}
 Then
\begin{eqnarray*}
\partial^*(S)&\ge&
\sum_{j=0}^{w-1}\left(\frac{w-j}{h-j}q_{h-j}(1-t_{h-j})+\frac{w-j}{h-j-1}q_{h-j}t_{h-j}\right)\\
&=&\sum_{j=0}^{w-1}\frac{w-j}{h-j}q_{h-j}+\sum_{j=0}^{w-1}\frac{w-j}{(h-j)(h-j-1)}q_{h-j}t_{h-j}\\
&\ge&\sum_{j=0}^{w-1}\frac{w-j}{h-j}q_{h-j}+\frac{1}{h^2}\sum_{j=0}^{w-1}(w-j)q_{h-j}t_{h-j}.
\end{eqnarray*}
If $\sum_{j=0}^{w-1}\frac{w-j}{h-j}q_{h-j}\ge (1-\delta)ks$ for some
$\delta>0$, then ~\eqn{eq:FlowCondition1} implies that
$$
\frac{1}{h^2}\sum_{j=0}^{w-1}(w-j)q_{h-j}t_{h-j}<\delta ks.
$$
Therefore $\eta(S)\le \sum_{j=0}^{w-1}q_{h-j}t_{h-j}<h^2\delta ks$.
This proves part (iv).\qed \ss

\section{Proof of Theorem~\ref{t:ForCore}}
\lab{probabilistic}

 Recall from
Section~\ref{mainResult} that $\mathcal{M}(n,{\bf m},k+1)$ is
$\mathcal{M}_{n,{\bf m}}$, which is a random multihypergraph with
given edge sizes,  restricted to multihypergraphs with minimum
degree at least $k+1$. In this section we prove the only remaining
theorem, Theorem~\ref{t:ForCore}. This theorem relates the
orientability of $\mathcal{M}(n,{\bf m},k+1)$ to its $w$-density.
  Recall that this probability space
was important because, by Proposition~\ref{uniformCore}, it gives
the distribution of the $(w,k+1)$-core $\coreH$ of
$\orH\in\mathcal{M}_{\bar n,\bar m,h}$ conditioned on the values of
$n$, the number of vertices and $m_{h-j}$, the number of hyperedges
of size $h-j$ for each $j$, in the core.

 It is  clear, that given values of
$n$ and ${\bf m}$, the probability space of random multihypergraphs
generated by $\mathcal{P}([n],{\bf M},{\bf 0},k+1)$, with
$|M_{h-j}|= (h-j)m_{h-j}$ ($h=0,\ldots,w-1$), is equivalent to
$\mathcal{M}(n,{\bf m},k+1)$. So we may, and do, make use of the
partition-allocation model for proving results about
$\mathcal{M}(n,{\bf m},k+1)$.

For the rest of the chapter, let $\eps>0$ and $k\ge 2$ be fixed.
Without loss of generality, we may assume that $\eps<\frac{1}{2}$
since $\eps$ may be taken arbitrarily small. By the hypothesis of
Theorem~\ref{t:ForCore}, we consider only ${\bf m}$ such that
$\sum_{j=0}^{w-1}(w-j)m_{h-j}\le kn-\eps
 n$. We may also assume that $\sum_{j=0}^{w-1}(w-j)m_{h-j}\ge kn-2\eps
 n$ since otherwise, by Theorem~\ref{t:CoreDensity}, we can simply add a set of random hyperedges so that the
 assumption holds. This is valid because   $(w,k)$-orientability is
 a decreasing property (i.e.\ it holds in all subgraphs of   $G$ whenever  $G$ has the property).
 Let
\begin{equation}
 D=\sum_{j=0}^{w-1}(h-j)m_{h-j},\ \  m=\sum_{j=0}^{w-1}m_{h-j},\ \   \mu=\frac{D}{n}.\lab{eq:D}
\end{equation}
Since
$$ D\cdot\frac{1}{h-w+1}\le
\sum_{j=0}^{w-1}(w-j)m_{h-j}\le D\cdot\frac{w}{h},\ \ \
m\le\sum_{j=0}^{w-1}(w-j)m_{h-j}\le wm,
$$
and
\begin{equation}
kn-2\eps n\le\sum_{j=0}^{w-1}(w-j)m_{h-j}\le kn-\eps n, \lab{eq:m}
\end{equation}
we have
\begin{equation}
 h(k-1)/w\le \mu=D/n\le (h-w+1)k,\ \ \ \frac{(k-1)n}{w}\le m\le \Big(k-\frac{1}{2}\Big)n. \lab{eq:mu-range}
\end{equation}

In the rest of the paper, whenever we refer to the probability space
$\mathcal{H}(n,{\bf m},k+1)$ or $\mathcal{M}(n,{\bf m},k+1)$, we
assume ${\bf m}$ satisfies~\eqn{eq:m}.

 Since $h$ and $w$ are given, we consider them as
absolute constants. Therefore, whenever we refer to $g=O(f)$, it
means that there exists a constant $C$ such that $g\le Cf$, where
$C$ can depend on $h$ and $w$. We also use notation
$g=O_{\gamma}(f)$, which means that there exists a constant $C$
depending on $\gamma$ only such that $g\le Cf$. The same convention
applies to $o(f)$, $\Omega(f)$, $\Theta(f)$ and $o_{\gamma}(f)$,
$\Omega_{\gamma}(f)$ and $\Theta_{\gamma}(f)$.

By Proposition~\ref{uniformCore}, conditioned on the values of $n$,
the number of vertices and $m_{h-j}$, the number of hyperedges of
size $h-j$, of the $(w,k+1)$-core $\coreH$ of
$\orH\in\mathcal{M}_{\bar n,\bar m,h}$, $\coreH$ is distributed as
$\mathcal{M}(n,{\bf m},k+1)$. Recall that $\eps>0$ and $k\ge 2$ are
fixed. Let ${\bf m}$ be an integer vector with the
constraint~\eqn{eq:m}. Given ${\bf m}$, let $D$, $\mu$ be as defined
in~\eqn{eq:D}. Let $\core2H$ be a random multihypergraph from the
probability space $\mathcal{M}(n,{\bf m},k+1)$.

We next sketch the proof of Theorem~\ref{t:ForCore}.  Let
$q_{h-j}(S)$ and $\eta(S)$ be defined as in~\eqn{eq:eta}. The
partition-allocation model  gives a good foundation for proving that
a.a.s.\ certain properties hold concerning the distribution of
vertex degrees and intersections of hyperedge sets with vertex sets.
Using this and various other probabilistic tools, we show that
\begin{description}

\item{(a)} the probability that $\core2H\in\mathcal{M}(n,{\bf m},k+1)$ has property $\mathcal{A}(\gamma)$ and contains some set $S$ with $|S|<\gamma n$ for which both Lemma~\ref{l:deterministic}(ii) and (iii) holds is
$o(1)$;

\item{(b)}  
there exists $\delta>0$, such that when $k$ is large enough, a.a.s.\
$\sum_{j=0}^{w-1}\frac{w-j}{h-j}q_{h-j}\ge (1-\delta)k|S|$,  and the
probability of $\core2H\in\mathcal{M}(n,{\bf m},k+1)$ containing
some set $S$ with $\gamma n\le |S|\le (1-\gamma)n$ and
$\eta(S)<h^2\delta k|S|$ is $o(1)$.

\end{description}
We also show the deterministic result that
\begin{description}
\item{(c)} no multihypergraph $\core2H$ with property $\mathcal{A}(\gamma)$ contains any sets $S$ with $|S|>(1-\gamma) n$ for which Lemma~\ref{l:deterministic}(i)
holds.
\end{description}
It follows that the probability that $\core2H$ has property
$\mathcal{A}(\gamma)$ and contains some set $S$ for which all parts
(i)--(iv) of Lemma~\ref{l:deterministic} hold is $o(1)$. Then by
Lemmas~\ref{l:expansion} and~\ref{l:deterministic},
$$
 \pr(\core2H\in
\mathcal{A}(\gamma)\wedge \core2H \ \mbox{is not}\
(w,k)\mbox{-orientable})=o(1).
$$

Finally, Lemma~\ref{uniformSimple} shows that the result applies to
random (simple) hypergraphs as well.

We start with a few concentration properties.  As discussed in
Section~\ref{core}, the degree sequence of
$\core2H\in\mathcal{M}(n,{\bf m},k+1)$ obeys the multinomial
distribution. The following lemma bounds the probability of rare
degree (sub)sequences where the degree distribution is independent
truncated Poisson. We will use this result to bound the probability
of rare degree sequences in $\mathcal{M}(n,{\bf m},k+1)$.
\begin{lemma}\lab{l:Concentration3} Let $s\ge w(n)$ for some $w(n)\rightarrow\infty$
as $n\rightarrow\infty$ and let $Y_1,\ldots, Y_s$ be independent copies of $Z$ defined
in~\eqn{eq:truncatedPoisson} with $\lambda$ satisfying $\lambda
f_k(\lambda)=\mu f_{k+1}(\lambda)$. Let  $0<\delta<1$ be any
constant. Then there exist $N>0$ and $0<\alpha<1$ both depending
only on $\delta$, such that provided $k>N$,
$$
\pr\left(\Big|\sum_{i=1}^s Y_i-\mu s\Big|\ge \delta\mu s\right)\le
\alpha^{\mu s}.
$$
\end{lemma}


\proof
Let $G(x)$ be the probability generating function of $Y_i$. Then
$$
G(x)=\sum_{j\ge k+1}
\pr(Z=j)x^j=\frac{e^{-\lambda}}{f_{k+1}(\lambda)}\left(e^{\lambda
x}-\sum_{j=0}^k\frac{(\lambda x)^j}{j!}\right)\le\frac{e^{\lambda
x-\lambda}}{f_{k+1}(\lambda)},
$$
for all $x\ge 0$. For any nonnegative integer $\ell$,
$$
\pr\left(\sum_{i=0}^s Y_i=l\right)\le \frac{G(x)^s}{x^{\ell}},\ \
\forall x\ge 0.
$$
Putting $x=\ell/s\lambda$ gives
\begin{equation}
\pr\left(\sum_{i=0}^s Y_i=\ell\right)\le \frac{e^{\ell-\lambda
s}}{(\ell/(\lambda
s))^{\ell}f_{k+1}(\lambda)^s}=\left(\frac{es\lambda}{\ell}\right)^{\ell}\left(\frac{e^{-\lambda}}{f_{k+1}(\lambda)}\right)^s.\lab{eq:1}
\end{equation}
It is easy to check that the right hand side of~\eqn{eq:1} is an
increasing function of $l$ when $l\le \lambda s$ and decreasing
function of $l$ when $l\ge \lambda s$. By
Proposition~\ref{p:lambda-mu}, there exists a constant $N_0$
depending only on $\delta$ such that provided $k>N_0$,
$(1-\delta)\mu<\lambda $. Thus, for any $\ell\le (1-\delta)\mu s$,
$$
\pr\left(\sum_{i=0}^s Y_i=\ell\right)\le
\left(\frac{es\lambda}{(1-\delta)\mu s}\right)^{(1-\delta)\mu
s}\left(\frac{e^{-\lambda}}{f_{k+1}(\lambda)}\right)^s,
$$
and so
 $$
\pr\left(\sum_{i=1}^s Y_i\le (1-\delta)\mu s\right)\le \mu
s\left(\frac{e\lambda}{(1-\delta)\mu }\right)^{(1-\delta)\mu
s}\left(\frac{e^{-\lambda}}{f_{k+1}(\lambda)}\right)^s.
$$
The expectation of $Y_1$ is  $\lambda
 f_k(\lambda)/f_{k+1}(\lambda)=\mu$. By Proposition~\ref{p:lambda-mu},
 we have $\mu\ge \lambda$ and $\mu-\lambda\rightarrow 0$ as $k\rightarrow \infty$.  Therefore
\begin{eqnarray*}
\pr\left(\sum_{i=1}^s Y_i\le (1-\delta)\mu s\right)&\le& \mu
s\left(\frac{\exp(\mu-\lambda-\delta\mu)}{(1-\delta)^{(1-\delta)\mu}f_{k+1}(\lambda)}\right)^s\\
&=&\mu
s\left(\frac{\exp(\mu-\lambda)}{f_{k+1}(\lambda)}\cdot\left(\frac{\exp(-\delta)}{(1-\delta)^{(1-\delta)}}\right)^{\mu}\right)^s.
\end{eqnarray*}
Since $0<\delta<1$,
$$
0<\frac{\exp(-\delta)}{(1-\delta)^{(1-\delta)}}<1.
$$
Since
$$
\exp(\mu-\lambda)\rightarrow 1,\ \ f_{k+1}(\lambda)\rightarrow 1,\ \
\mbox{as}\ k\rightarrow \infty
$$
by Proposition~\ref{p:lambda-mu}, there exists $N_1>0$ and
$0<\alpha_1<1$, both depending only on $\delta$, such that provided
$k>N_1$,
$$
\pr\left(\sum_{i=1}^s Y_i\le (1-\delta)\mu s\right)\le \alpha_1^{\mu
s}.
$$
Now we bound the upper tail of $\sum_{i=1}^{s}Y_i$. Let
$j=1,2,\ldots$. For any $\ell$ satisfying $(1+j)\mu s\le \ell<
(2+j)\mu s$, as with the lower tail bound,
\begin{equation}
\pr\left(\sum_{i=0}^s Y_i=\ell\right)\le
\left(\frac{es\lambda}{(1+j)\mu s}\right)^{(1+j)\mu
s}\left(\frac{e^{-\lambda}}{f_{k+1}(\lambda)}\right)^s=
\left(\frac{\exp(\mu-\lambda)}{f_{k+1}(\lambda)}\cdot\left(\frac{e^j}{(1+j)^{(1+j)}}\right)^{\mu}\right)^s,\nonumber
\end{equation}
and so
\begin{equation}
\pr\left((1+j)\mu s\le\sum_{i=1}^s Y_i< (2+j)\mu s\right)\le \mu
s\left(\frac{\exp(\mu-\lambda)}{f_{k+1}(\lambda)}\cdot\left(\frac{e^j}{(1+j)^{(1+j)}}\right)^{\mu}\right)^s.\lab{eq:2}
\end{equation}
 Similarly we have
\begin{equation}
\pr\left((1+\delta)\mu s\le\sum_{i=1}^s Y_i< 2\mu s\right)\le \mu
s\left(\frac{\exp(\mu-\lambda)}{f_{k+1}(\lambda)}\cdot\left(\frac{e^{\delta}}{(1+\delta)^{(1+\delta)}}\right)^{\mu}\right)^s.\lab{eq:3}
\end{equation}
Since $ 0<e^{\delta}/(1+\delta)^{(1+\delta)}<1$ for any $\delta>0$,
we may bound the right side of~\eqn{eq:2} and~\eqn{eq:3} by
$\alpha_2^{\mu s}$ where $0<\alpha_2<1$ is some constant depending
only on $\delta$. Also, since $e/(1+j)<1$ for all $j\ge 2$, the
right side of~\eqn{eq:2} is at most $\exp(-\Omega(\mu)js)$ for $j\ge
2$ provided $k$ is large enough.
 Hence there exists $N_2>0$ and
$0<\alpha_3<1$ depending only on $\delta$, such that provided
$k>N_2$,
\begin{eqnarray*}
\pr\left(\sum_{i=1}^s Y_i\ge (1+\delta)\mu s\right)\le\alpha_3^{\mu
s}.
\end{eqnarray*}
The lemma follows by choosing $\alpha=\max\{\alpha_1,\alpha_3\}$ and
$N=\max\{N_1,N_2\}$.
 \qed
\begin{lemma}\lab{l:prob}
Let $k\ge -1$ be an integer. Drop $D$ balls independently at random
into $n$ bins. Let $\mu=D/n$ and let $\lambda$ be defined as $
\lambda f_k(\lambda)=\mu f_{k+1}(\lambda)$. Assume
$D-(k+1)n\rightarrow\infty$ as $n\rightarrow\infty$. Then the
probability that each bin contains at least $k+1$ balls is
$\Omega(f_{k+1}(\lambda)^n)$.
\end{lemma}


\proof Let ${\bf d}$ denote $(d_1,\ldots,d_n)$.  Let
$\mathscr{D}=\{{\bf d}:d_i\ge k+1\ \forall i\in[n],\ \sum_{i=1}^n
d_i=D\}$. Let $\pr(B)$ denote the probability that each bin contains
at least $k+1$ balls. Then
$$
\pr(B)=\sum_{{\bf
d}\in\mathscr{D}}\binom{D}{d_1,\ldots,d_n}\Big/n^D=\frac{D!}{n^D}\sum_{{\bf
d}\in\mathscr{D}}\prod_{i=1}^n\frac{1}{d_i!}.
$$
Let $Y_1,\ldots, Y_n$ be $n$ independent truncated Poisson variables
which are copies of $Z_{(\ge k+1)}$  as defined
in~\eqn{eq:truncatedPoisson} with parameter $\lambda$ satisfying $
\lambda f_k(\lambda)=\mu f_{k+1}(\lambda)$. Then
$$
\pr\left(\sum_{i=1}^n Y_i=D\right)=\sum_{{\bf
d}\in\mathscr{D}}\prod_{i=1}^n\frac{e^{-\lambda}\lambda^{d_i}}{f_{k+1}(\lambda)d_i!}=\frac{e^{-\lambda
n}\lambda^D}{f_{k+1}(\lambda)^n}\sum_{{\bf
d}\in\mathscr{D}}\prod_{i=1}^n\frac{1}{d_i!}.
$$
Since $D-(k+1)n\rightarrow\infty$ as $n\rightarrow\infty$,
$\pr\left(\sum_{i=1}^n Y_i=D\right)=\Omega(D^{-1/2})$
(see~\cite[Theorem 4(a)]{PW} for a short proof),
$$
\sum_{{\bf
d}\in\mathscr{D}}\prod_{i=1}^n\frac{1}{d_i!}=\Omega\left(\frac{e^{\lambda
n}f_{k+1}(\lambda)^n}{\lambda^D D^{1/2}}\right).
$$
So, using Stirling's formula,
\begin{eqnarray}
\pr(B)&=&\Omega\left(\frac{D!}{n^D}\cdot\frac{e^{\lambda
n}f_{k+1}(\lambda)^n}{\lambda^D D^{1/2}}\right)=\Omega\left(\sqrt{
D}\left(\frac{D}{en}\right)^D\cdot\frac{e^{\lambda
n}f_{k+1}(\lambda)^n}{\lambda^D D^{1/2}}\right)\lab{eq:prob2}\\
&=&\Omega\left(\left(\frac{\mu}{\lambda}e^{\lambda/\mu-1}\right)^{\mu
n}f_{k+1}(\lambda)^n\right). \nonumber
\end{eqnarray}
Since $(\mu/\lambda)\cdot e^{\lambda/\mu-1}\ge 1$, $
\pr(B)=\Omega(f_{k+1}(\lambda)^n) $.\qed
\begin{cor}\lab{c:prob}
Let $k\ge -1$ be an integer. Let $\mathscr{D}=\{{\bf d}:d_i\ge k+1,
\forall i\in[n], \sum_{i=1}^n d_i=D\}$ and let $A_n$ be any subset
of $\mathscr{D}$. Let $\mu=D/n$. Let $\pr(A_n)$ denote the
probability that the degree sequence ${\bf d}$ of
$\core2H\in\mathcal{M}(n,{\bf m},k+1)$ is in $A_n$ and let
$\pr_{TP}(A_n)$ be the probability that $(Y_1,\ldots,Y_n)\in A_n$
where $Y_i$ are independent copies of the random variable $Z_{(\ge
k+1)}$ as defined in~\eqn{eq:truncatedPoisson} with the parameter
$\lambda$ satisfying $\lambda f_k(\lambda)=\mu f_{k+1}(\lambda)$.
Assume $D-(k+1)n\rightarrow\infty$ as $n\rightarrow\infty$. Then
$$
\pr(A_n)=O\left(\sqrt{D}\right)\pr_{TP}(A_n).
$$

\end{cor}


\proof Let $A_n$ be any subset of $\mathscr{D}$ and let $\pr(B)$
denote the probability that each bin contains at least $k+1$ balls
by dropping $D$ balls independently and randomly into $n$ bins.
Consider the partition-allocation model that generates
$\mathcal{M}(n,{\bf m},k+1)$, which allocates the partitioned $D$
balls randomly into $n$ bins with the restriction that each bin
contains at least $k+1$ balls. Then
\begin{eqnarray*}
\pr(A_n)=\sum_{{\bf d}\in A_n}\frac{1}{\pr(B)}\cdot
\binom{D}{d_1,\ldots,d_n}\Big/ n^D=\frac{D!}{n^D \pr(B)}\sum_{{\bf
d}\in A_n}\prod_{i=1}^n\frac{1}{d_i!},
\end{eqnarray*}
and
\begin{eqnarray*}
\pr_{TP}(A_n)=\sum_{{\bf d}\in
A_n}\prod_{i=1}^n\frac{e^{-\lambda}\lambda^{d_i}}{f_{k+1}(\lambda)d_i!}=\frac{e^{-\lambda
n}\lambda^D}{f_{k+1}(\lambda)^n}\sum_{{\bf d}\in
A_n}\prod_{i=1}^n\frac{1}{d_i!}.
\end{eqnarray*}
Therefore
\begin{eqnarray*}
&&\pr(A_n)=\frac{D!e^{\lambda n}f_{k+1}^n}{n^D
\pr(B)\lambda^D}\pr_{TP}(A_n)=O\left(\sqrt{D}\right)\pr_{TP}(A_n),
\end{eqnarray*}
since $
\pr(B)=\Omega\left(\left(\frac{\mu}{\lambda}e^{\lambda/\mu-1}\right)^{\mu
n}f_{k+1}(\lambda)^n\right)$  by Lemma~\ref{l:prob}~\eqn{eq:prob2}.
\qed \ss

A significant difficulty in this work is to ensure that various
constants do not depend on the choice of $\eps$. In particular, we
emphasize that the constants such as $\alpha$ and $N$ in the
following results do not depend on $\eps$.

The next is a corollary of Lemma~\ref{l:Concentration3} and
Corollary~\ref{c:prob}.
\begin{cor}\lab{c:Concentration3}  Let $\mu$ be defined as in~\eqn{eq:D}. Let $0<\delta<1$ be any constant.
Then there exist two constants $N>0$ and $0<\alpha<1$, both
depending only on $\delta$, such that, provided $k>N$, for any
vertex set $S\subset V(\core2H)$  with $|S|\ge \log^2 n$,
$$
\pr(|d(S)-\mu |S||\ge \delta \mu |S|)\le\alpha^{\mu|S|}.
$$
\end{cor}


\proof  Let $Y_1,\ldots,Y_n$ be independent copies of the truncated
Poisson random variable $Z$ as defined in~\eqn{eq:truncatedPoisson}.
Let $S\subset V(\core2H)$ and let $s=|S|$. Then by
Lemma~\ref{l:Concentration3}, there exist $N>0$ and
$0<\hat\alpha<1$, both depending only on $\delta$, such that
provided $k>N$,
\begin{equation}
\pr\left(\Big|\sum_{i\in S} Y_i-\mu s\Big|\ge \delta \mu s
\right)\le \hat\alpha^{\mu s},\nonumber
\end{equation}
 By Corollary~\ref{c:prob},
\begin{eqnarray*}
\pr(|d(S)-\mu s|\ge \delta \mu s) &\le&O(D^{1/2})\hat\alpha^{\mu
s}=\left(\exp\left(\frac{\ln \Theta(\sqrt{\mu n})}{\mu
s}\right)\hat\alpha\right)^{\mu s}.
\end{eqnarray*}
Since  $s\ge \log^2 n$ and so
$$
\frac{\ln \Theta(\sqrt{\mu n})}{\mu s}\rightarrow 0, \ \ \mbox{as}\
n\rightarrow\infty.
$$
Let $\alpha=1/2+\hat\alpha/2$. Then $0<\hat\alpha<\alpha<1$ and
$\alpha$ depends only on $\delta$. Then provided $k>N$, $
\pr(|d(S)-\mu s|\ge \delta \mu s)\le \alpha^{\mu s} $. \qed \ss

The following corollary shows that $d(S)$ is very concentrated when
$S$ is not too small.
\begin{cor}\lab{c:Concentration} Let $\delta>0$ and $0<\gamma<1$ be arbitrary constants. Then there exists a constant $N$ depending only on
$\delta$ and $\gamma$, such that provided $k>N$,
$$
\pr(\exists S\subset V(\core2H), s\ge \gamma n, |d(S)-\mu
s|\ge\delta\mu s)=o(1).
$$

\end{cor}

\proof For any $S\subset V(\core2H)$, let $s=|S|$.  By
Corollary~\ref{c:Concentration3}, there exists $N_1>0$ and
$0<\alpha<1$, both depending only on $\delta$, such that provided
$k>N_1$, for any $S\subset V(\core2H)$,
$$
\pr(|d(S)-\mu s|\ge \delta \mu s)\le \alpha^{\mu s}.
$$
Let $N_2$ be the smallest integer such that
$e\alpha^{N_2}/\gamma<1/2$. Let $N=\max\{N_1,N_2\}$. Then $N$
depends only on $\delta$ and $\gamma$. For all $\mu>N$,

\begin{eqnarray*}
&&\pr(\exists S\subset [n], s\ge \gamma n, |d(S)-\mu s|\ge\delta\mu
s)\le \sum_{\gamma n\le s\le n}\binom{n}{s}\alpha^{\mu
s}\le\sum_{\gamma n\le s\le n}\left(\frac{en}{s}\cdot\alpha^{\mu
}\right)^s\\
&&\hspace{0.4cm}\le\sum_{\gamma n\le s\le
n}\left(\frac{e}{\gamma}\cdot\alpha^{\mu
}\right)^s=O\left(2^{-\gamma n}\right)=o(1).\qed
\end{eqnarray*}

The following lemma will be used later to prove that a.a.s.\
$\sum_{j=0}^{w-1}(w-j)q_{h-j}/(h-j)\ge (1-\delta)ks$ provided $k$ is
large enough.
\begin{lemma}\lab{l:Concentration4}
Let $\mathcal{C}=\{c_0,\ldots,c_{w-1}\}$ be a set of colours.
Suppose that  $D$ balls are each coloured with some colour in
$\mathcal{C}$, and let  $p_j$ denote the proportion of balls that
are coloured $c_j$ ($0\le j\le w-1$). Randomly choose a subset of
$q$ of the balls.
 Let $q_{j}$ be the
number of balls chosen that are coloured with $c_j$. Then for any
$0\le j\le w-1$ and $0<\delta < 1$,
$$
\pr(|q_{j}-p_jq|\ge \delta p_jq)\le \exp(-\Omega(\delta^2 p_jq)).
$$
\end{lemma}

\proof For any $0\le j\le w-1$ any $\ell>0$,
$$
\pr(q_j=\ell)=\binom{p_jD}{\ell}\binom{D-p_jD}{q-\ell}\Big/\binom{D}{q}.
$$
Let $p_{\ell}$ denote $\pr(q_j=\ell)$. Put $\ell_0=p_j q$,
$\ell_1=(1-\delta/2)p_j q$ and $\ell_2=(1-\delta)p_j q$. Then for
any $\ell\le \ell_1$,
$$
\frac{p_{\ell-1}}{p_{\ell}}=\frac{\ell(D(1-p_j)-q+\ell)}{(p_jD-\ell+1)(q-\ell
+1)}\le\frac{\ell_1(D(1-p_j)-q+\ell_0)}{(p_jD-\ell_0)(q-\ell_0
)}=1-\frac{\delta}{2}.
$$
Then
$$
p_{\ell_2}\le (1-\delta/2)^{\delta p_j
q/2}p_{\ell_1}\le\exp\left(\frac{\delta p_j q}{2}\ln
\Big(1-\frac{\delta}{2}\Big)\right)\le\exp(-\delta^2 p_j q/4).
$$
So
$$
\pr(q_{j}\le (1-\delta)p_jq)=\sum_{\ell\le\ell_2}p_{\ell}\le
\frac{1}{\delta}p_{\ell_2}\le \exp(-\Omega(\delta^2 p_jq)).
$$
Similarly we can bound the upper tail and then
Lemma~\ref{l:Concentration4} follows. \qed
\begin{lemma}\lab{l:Concentration2} Let $0<\delta<1$ and $0<\gamma<1$ be two arbitrary
constants. Given $S\subset V(\core2H)$, let $q_{h-j}=q_{h-j}(S)$ be
as defined in~\eqn{eq:eta}. Then there exists $N>0$ depending only
on $\delta$ and $\gamma$ such that for all $k>N$,
$$
\pr\left(\exists S\subset V(\core2H), |S|\ge\gamma n,
\sum_{j=0}^{w-1}\frac{w-j}{h-j}q_{h-j}< (1-\delta)k|S|\right)=o(1).
$$

\end{lemma}


\proof For any $0\le j\le w-1$, let $p_j$ denoted $(h-j)m_{h-j}/D$.
Let $J:=\{j: p_j>\delta/8w\}$. We first show that given $S\subset
V(\core2H)$ with $|S|\ge\gamma n$, if
\begin{equation}
\sum_{j=0}^{w-1}\frac{w-j}{h-j}q_{h-j}< (1-\delta)k|S|,\lab{eq:comb}
\end{equation}
then there exists $j\in J$ such that $q_{h-j}(S)\le(1-\delta/8)p_j
d(S)$. Assume there is no such $j$ by contradiction. Then
\begin{eqnarray}
&&\hspace{-0.5cm}\sum_{j=0}^{w-1}\frac{w-j}{h-j}q_{h-j}(S)\ge\sum_{j\in
J}\frac{w-j}{h-j}q_{h-j}(S) >(1-\delta/8)d(S)\sum_{j\in
J}\frac{w-j}{h-j}p_j\nonumber\\
&&\hspace{-0.5cm}\hspace{0.4cm}=(1-\delta/8)d(S)\left(\sum_{j=0}^{w-1}\frac{w-j}{h-j}p_j-\sum_{j\notin
J}\frac{w-j}{h-j}p_j\right)\ge(1-\delta/8)d(S)\left(\sum_{j=0}^{w-1}\frac{w-j}{h-j}p_j-\frac{w}{h}\frac{\delta}{8w}\right)\nonumber\\
&&\hspace{-0.5cm}\hspace{0.4cm}\ge(1-\delta/8)d(S)\sum_{j=0}^{w-1}\frac{w-j}{h-j}p_j(1-\delta/8)\ge(1-\delta/4)d(S)
\sum_{j=0}^{w-1}\frac{w-j}{h-j}p_j.\lab{eq:11}
\end{eqnarray}
Let $s=|S|$ and let $r=s/n$. Then by
Corollary~\ref{c:Concentration}, there exists $N_2>0$ depending on
$\delta$ and $\gamma$ only, such that a.a.s. $d(S)\ge
(1-\delta/4)Dr$ whenever $k>N_2$. Therefore, combining
with~\eqn{eq:11}, we get a.a.s.\ provided $k>\max\{N_1,N_2\}$,
\begin{eqnarray*}
\sum_{j=0}^{w-1}\frac{w-j}{h-j}q_{h-j}(S)>(1-\delta/2)Dr\sum_{j=0}^{w-1}\frac{w-j}{h-j}p_j\ge
(1-\delta/2)r(kn-2\eps n)=(1-\delta/2)(k-2\eps)s.
\end{eqnarray*}
For any $k>2/\delta\ge4\eps/\delta$, we have
$(1-\delta/2)(k-2\eps)s>(1-\delta)ks$. Take
$N=\max\{N_1,N_2,2/\delta\}$. Then for any $k>N$, we have a.a.s.\
$$
\sum_{j=0}^{w-1}\frac{w-j}{h-j}q_{h-j}(S)>(1-\delta)ks,
$$
which contradicts~\eqn{eq:comb}. It follows that there exists $j\in
J$ such that $q_{h-j}(S)\le(1-\delta/8)p_j d(S)$.

Consider the partition-allocation model that generates
$\mathcal{P}([n],{\bf M},{\bf 0},k+1)$. Let
$\mathcal{C}=\{c_0,\ldots,c_{w-1}\}$ be a set of colours. For balls
partitioned into parts that are of size $h-j$ for some $0\le j\le
w-1$, colour them with $c_j$. Then the $w$ colours are distributed
u.a.r.\ among the $D$ balls. By Lemma~\ref{l:Concentration4}, for
any $S\subset V(\core2H)$,
$$
\pr\big(q_{h-j}(S)\le (1-\delta/8) p_jd(S)\big)\le
\exp\big(-\Omega(\delta^2 p_jd(S))\big).
$$
 Then
there exists a constant $N_1$ depending only on $\delta$ and
$\gamma$ such that,
\begin{eqnarray*}
&&\pr\big(\exists S,j\in J, s\ge \gamma n, q_{h-j}(S)\le
(1-\delta/8)p_jd(S)\big)\\
&&\hspace{0.6cm}\le w2^n \exp\big(-\Omega(\delta^3 d(S))\big)\le
w\left(2 \exp(-\Omega(\delta^3 \gamma k))\right)^n=o(1).
\end{eqnarray*}
Note that the inequality holds because $|J|\le w$, the number of
sets $S$ with $|S|\ge \gamma n$ is at most $2^n$,  $\delta/8w\le
p_j<1$ for all $j\in J$ and $d(S)\ge (k+1)|S|>k\gamma n$. It follows
that a.a.s.\ there exists no set $S$ with $|S|\ge \gamma n$ for
which there exists $j\in J$ such that $q_{h-j}(S)\le(1-\delta/8)p_j
d(S)$. Lemma~\ref{l:Concentration2} then follows.\qed \ss

Recall that $\rho(S)$ is the number of hyperedges partially
contained in $S$ and $\nu(S)$ is the number of hyperedges
intersecting $S$ by the definition above
Lemma~\ref{l:deterministic}.
\begin{lemma}\lab{l2:SmallSet}
Let  $\delta>0$ be any  constant and let $\mu=\mu(\core2H)=D/n$ as
defined in~\eqn{eq:D}. Then there exists a constant $N>0$  depending
only on $\delta$ such that, provided $k>N$, a.a.s.\ there exists no
$S\subset V(\core2H)$ for which  $\log^2 n\le|S|\le n$, $d(S)<
(1-\delta)\mu |S|$, and $\nu(S)<k|S|$.
\end{lemma}

This lemma will be proved after the proof of
Theorem~\ref{t:ForCore}.\ss

 \no {\bf Proof of Theorem~\ref{t:ForCore}.\ }
By Lemma~\ref{l:expansion} and~\ref{l:deterministic}, it is enough
to show that the expected number of sets $S$ contained in a
hypergraph $\core2H\in\mathcal{M}(n,{\bf m},k+1)$ with property
$\mathcal{A}(\gamma)$ for which all of Lemma~\ref{l:deterministic}
(i)--(iv) are satisfied is $o(1)$. We call a set $S\subset
V(\core2H)$ is interesting if it lies in a hypergraph $\core2H$ with
property $\mathcal{A}(\gamma)$. Let $X$ be the number of interesting
sets $S\subset V(\core2H)$ such that~\eqn{eq:FlowCondition1} holds.
Similarly, let $X_{< a}$ (or $X_{> b}$ or $X_{[a,b]})$ for any
$0<a<b<n$ denote the number of interesting $S\subset [n]$ such that
~\eqn{eq:FlowCondition1} holds and $|S|< a$ (or $|S|> b$ or $a\le
|S|\le b)$ respectively. For any set $S$ under discussion, let $s$
denote $|S|$ and $\bar s$ denote $|\barS|$.

{\em Case 1:} $s<\eps n/k$. By theorem's hypothesis
$$
\left(\sum_{j=0}^{w-1}(w-j)m_{h-j}\right)-kn<-\eps n,
$$
any $S$ satisfying~\eqn{eq:FlowCondition1} must satisfy
\begin{equation}
\partial^*(S)<ks-\eps n.\lab{eq:FlowCondition}
\end{equation}
When $s<\eps n/k$, $ks-\eps n<0$. However $\partial^*(S)\ge 0$ as
observed below~\eqn{eq:expansion}. Hence~\eqn{eq:FlowCondition}
cannot hold. Thus $X_{< \eps n/k}=0$.

{\em Case 2:} $s>(1-\gamma)n$. part (i) of
Lemma~\ref{l:deterministic} says that \eqn{eq:FlowCondition1} holds
only if the number of hyperedges partially contained in $\barS$ is
at least $k\bar s/w$. But $X$ counts only interesting sets, i.e.\
sets that lie in a hypergraph with property $\mathcal{A}(\gamma)$.
By the definition of property $\mathcal{A}(\gamma)$, there are no
such interesting sets and so $X_{\ge(1-\gamma) n}=0$.


{\em Case 3: } $\eps n/k\le s<\gamma n$. Let $\delta_1=(h-w)/2h$.
By Lemma~\ref{l2:SmallSet}, there exists $N_1>0$ such that provided
$k>N_1$, the expected number of $S$ with $d(S)<(1-\delta_1)\mu s$
for which Lemma~\ref{l:deterministic} (ii) is satisfied and $\eps
n/k\le s\le n$ is $o(1)$. We now show that there exists no
interesting sets $S\subset V(\core2H)$ with $|S|<\gamma n$ for which
Lemma~\ref{l:deterministic} (iii) holds and $d(S)\ge(1-\delta_1)\mu
s$.
 If
$d(S)\ge(1-\delta_1)\mu s$, $d(S)\ge \frac{h+w}{2w}ks$ provided
$k\ge h+w$ since $\mu\ge h(k-1)/w$ by~\eqn{eq:mu-range}. Then it
follows that
$$
\frac{d(S)-ks}{h-w}\ge\frac{ks}{2w}.
$$
Lemma~\ref{l:deterministic} (iii) implies
that~\eqn{eq:FlowCondition1} holds only if the number of hyperedges
partially contained in $S$ is at least $ks/2w$. By the definition of
property $\mathcal{A}(\gamma)$, there is no such interesting sets
$S$ when $s<\gamma n$. So provided $k> \max\{N_1,h+w\}$, a.a.s.\
there exists no interesting sets $S$, with $s<\gamma n$ for which
both Lemma~\ref{l:deterministic} (ii) and (iii) hold. Then $\ex(X_{<
\gamma n})=o(1)$.

 Note that $k|S|/2w$ in the definition of property
$\mathcal{A}(\gamma)$ can be modified to be $Ck|S|$ for any positive
constant $C$, and it can be checked straightforwardly that there
exists a constant $\gamma$ depending on $C$ only, such that
Corollary~\ref{c:SmallSet2} holds. Therefore, any $0<\delta_1<1-w/h$
would work here by choosing some appropriate $C$ to modify the
definition of property $\mathcal{A}(\gamma)$.

{\em Case 4:} $\gamma n\le s\le (1-\gamma)n$. Let $0<\delta_2<1$ be
chosen later.
By Lemma~\ref{l:Concentration2}, there exists $N_2>0$ depending only
on $\delta_2$ such that provided $k\ge N_2$, a.a.s.
$$
\sum_{j=0}^{w-1}\frac{w-j}{h-j}q_{h-j}\ge (1-\delta_2)ks \ \ \
\mbox{for all}\ S \ \mbox{with}\ \gamma n\le |S|\le (1-\gamma)n.
$$
 For any $S\subset V(\core2H)$, let $\eta=\eta(S)$ be
as defined in~\eqn{eq:eta}. Then by Lemma~\ref{l:deterministic}, to
show $\ex(X_{[\gamma n, (1-\gamma)n]})=o(1)$, it is enough to show
that the expected number of sets $S$ with $\gamma n\le s\le
(1-\gamma) n$ for which $\eta(S)$ is at most $h^2\delta_2 ks$, is
$o(1)$. Consider the probability space  $\mathcal{M}(n,{\bf m},0)$,
which is generated by placing each hyperedge uniformly and randomly
on the $n$ vertices. Let $B$ be the event that all bins contain at
least $k+1$ balls. Then $\mathcal{M}(n,{\bf m},k+1)$ equals
$\mathcal{M}(n,{\bf m},0)$ conditioned on the event $B$. By
Lemma~\ref{l:prob} $\pr(B)=\Omega(f_{k+1}(\lambda)^n)$ where
$\lambda f_{k}(\lambda)=\mu f_{k+1}(\lambda)$. Given any set $S$,
let $r=s/n$. For any hyperedge of size $h-j$, the probability for it
to intersect both $S$ and $\barS$ is
$p_{j,r}=1-r^{h-j}-(1-r)^{h-j}$. Then $p_{j,r}\ge
1-\gamma^{h-j}-(1-\gamma)^{h-j}\ge
1-\gamma^{h-w+1}-(1-\gamma)^{h-w+1}$ for any set $S$ and any $0\le
j\le w-1$. Recall from~\eqn{eq:D} that $m$ is the total number of
hyperedges in $\core2H$. Then $\ex
\eta(S)=\sum_{j=0}^{w-1}p_{j,r}m_{h-j}\ge
m(1-\gamma^{h-w+1}-(1-\gamma)^{h-w+1})$  for any given $S$. Since
 $m\ge (k-1)n/w$ by~\eqn{eq:mu-range},
$$
\ex \eta(S)\ge
(1-\gamma^{h-w+1}-(1-\gamma)^{h-w+1})(k-1)n/w=\Theta_{\gamma}(k)n, \
 \ \mbox{for any}\ S\ \mbox{with}\ \gamma n\le|S|\le (1-\gamma)n.
$$
Choose
$$
\delta_2=\frac{1-\gamma^{h-w+1}-(1-\gamma)^{h-w+1}}{4wh^2(1-\gamma)}.
$$
Then $\delta_2$ depends only on $\gamma$ and so $N_2$  also depends
only on $\gamma$. By the Chernoff bound~\cite{C3},
\begin{eqnarray*}
\pr(\eta(S)<h^2\delta_2 ks)\le \pr(\eta(S)<h^2\delta_2
k(1-\gamma)n)\le \pr\left(\eta(S)<\frac{1}{2}\ex \eta(S)\right)\le
\exp(-\ex\eta(S)/16).
\end{eqnarray*}
Note that the second inequality holds because of the choice of
$\delta_2$. So there exists some constant $C>0$ s.t.
\begin{eqnarray*}
\pr(\eta(S)<h^2\delta_2 ks\mid B)\le
C\exp\left(-\ex\eta(S)/16\right)f_{k+1}(\lambda)^{-n}
=C\left(\exp\left(-\frac{\ex\eta(S)}{16n}-\ln
f_{k+1}(\lambda)\right)\right)^n.
\end{eqnarray*}
The number of sets $S$ with $\gamma n\le |S|\le (1-\gamma)n$ is at
most $2^n$. So the expected number of sets $S$ with $\gamma n\le
s\le (1-\gamma) n$ and $\eta(S)<h^2\delta_2 ks$ in
$\mathcal{M}(n,{\bf m},k+1)$ is at most
$$
C\left(2\exp\left(-\frac{\ex\eta(S)}{16n}-\ln
f_{k+1}(\lambda)\right)\right)^n.
$$
Clearly $f_{k+1}(\lambda)\rightarrow 1$ as $k\rightarrow\infty$ and
$\ex\eta(S)=\Theta_{\gamma}(k)n$ as observed before. Then there
exists a constant $N_3>0$ depending only on $\gamma$ such that
provided $k>N_3$,
$$
2\exp\left(-\frac{\ex\eta(S)}{16n}-\ln f_{k+1}(\lambda)\right)<1.
$$
Then provided $k\ge \max\{N_2,N_3\}$,
$$
\ex(X_{[\gamma n,(1-\gamma)n]})=o(1).
$$

Combining all cases, let $N=\max\{N_1,N_2,N_3,h+w\}$. Then $N$
depends only on $\gamma$. We have shown that provided $k> N$,
 $ \ex X=o(1) $. Then
Theorem~\ref{t:ForCore} follows. \qed \ss

\no {\bf Proof of Lemma~\ref{l2:SmallSet}.\ }   The idea of the
proof is as follows. When  $S$ is big, by
Corollary~\ref{c:Concentration} there are no such sets with $d(S)<
(1-\delta)\mu |S|$. We will see later that $\nu(S)<k|S|$ requires a
lot of hyperedges partially contained in $S$, which is unlikely to
happen when  $S$ is small enough.

 Let $\core2H\in\mathcal{M}(n,{\bf m},k+1)$. Let
$D=\sum_{j=0}^{w-1}(h-j)m_{h-j}$ and $\mu=D/n$ as defined
in~\eqn{eq:D}. For any $S$, let $\rho(S,i)$ denotes the number of
hyperedges with exactly $i$ vertices contained in $S$. Then
$\nu(S)<ks$ if and only if $\sum_{i=2}^{h}(i-1)\rho(S,i)>d(S)-ks$.
By Corollary~\ref{c:Concentration}, there exists $N_1>0$ depending
only on $\delta$ such that provided $k>N_1$, a.a.s.\ there is no $S$
such that $s>n/h$ and $d(S)<(1-\delta)\mu s$. So we only need to
consider sets $S$ with $|S|\le n/h$.  We call a vertex set $S\in
\core2H$ {\em bad} if $\log^2 n\le |S|\le n/h$, $d(S)<(1-\delta)\mu
|S|$ and $\sum_{i=2}^h(i-1)\rho(S,i)> d(S)-k|S|$. Let $s$ denote
$|S|$.

  For any given $S$, let $p(S)$
denote the probability of $S$ being bad. By
Corollary~\ref{c:Concentration3}, there exists $N_2>0$ and
$0<\alpha<1$, both depending only on $\delta$, such that provided
$k>N_2$, the probability that $d(S)<(1-\delta)\mu s$ is at most
$\alpha^{\mu s}$. Let $p(q,t)$ be the probability that
 that $\sum_{i=2}^h(i-1)\rho(S,i)$ is at least $t$ conditional on $d(S)=q$.
Then
\begin{equation}
p(S)=\sum_{(k+1)s\le q\le (1-\delta)\mu s}
p(q,q-ks)\pr(d(S)=q).\lab{eq:p(S)}
\end{equation}
For the small value of $q$ (or $s$), we need the following claim, to
be proved later.
\begin{claim}\lab{c:p(q)} If $q<D/h$, then
$$
p(q,t)\le \left(\exp\left(\frac{h\ln
t}{t}\right)\frac{eh(h-1)^2q^2}{4tD}\right)^t.
$$
In particular, if $t\rightarrow \infty$ as $n\rightarrow \infty$,
then
$$
p(q,t)\le \left(\frac{eh^3q^2}{4tD}\right)^t.
$$
\end{claim}
{\em Case 1: }   $s<2n/eh^3(k+1)$.
 Since $ (k+1)s\le q\le
(1-\delta)\mu s<D/h$, we have
$$
\frac{q}{q-ks}\le\frac{(k+1)s }{(k+1)s-ks}= k+1, \ \ \ \
\frac{q}{D}\le\frac{(1-\delta)\mu s}{\mu n}<\frac{s}{n},\ \ \
q-ks\ge s\ge\log^2 n.
$$
So $q-ks\rightarrow \infty$ as $n\rightarrow \infty$.
By~\eqn{eq:p(S)} and the particular case of Claim~\ref{c:p(q)}, we
have
\begin{eqnarray*}
p(S)&\le& \sum_{(k+1)s\le q\le (1-\delta)\mu s}
\left(\frac{eh^3(k+1)s}{4n}\right)^{q-ks}\pr(d(S)=q)\\
&\le& \left(\frac{eh^3(k+1)s}{4n}\right)^s\pr\Big((k+1)s\le d(S)\le
(1-\delta)\mu s\Big)\le
\left(\frac{eh^3(k+1)s}{4n}\right)^s\alpha^{\mu s}.
\end{eqnarray*}
Note that the second inequality above holds because $q-ks\ge s$ and
$0<eh^3(k+1)s/4n<1$ since $s<2n/eh^3(k+1)$.

Then the expected number of bad sets $S$ with $|S|=s$, for any fixed
$\log^2 n\le s<2n/eh^3(k+1)$, is at most
\begin{eqnarray*}
\binom{n}{s}\left(\frac{eh^3(k+1)s}{4n}\right)^s\alpha^{\mu s}\le
\left(\frac{en}{s}\cdot\alpha^{\mu}\cdot\frac{eh^3(k+1)s}{4n}\right)^s=\left(e^2h^3
(k+1)\alpha^{\mu}/4\right)^s.
\end{eqnarray*}
Since $\mu\ge h(k-1)/w$ by~\eqn{eq:mu-range}, this is at most
$\exp(-s)$ provided $k\ge N_3$ for some $N_3>0$ depending only on
$\alpha$.

{\em Case 2: } $s\ge 2n/eh^3(k+1)$. Take $p(q,q-ks)\le 1$ since
$p(q,q-ks)$ is a probability. So the expected number of bad sets $S$
with $|S|=s$, for any fixed $2n/eh^3(k+1)\le s\le (1-\delta)\mu s$,
is at most
$$
\binom{n}{s}\cdot \alpha^{\mu
s}=\left(\frac{en}{s}\alpha^{\mu}\right)^s\le\left(e^2h^3(k+1)\alpha^{\mu}/2\right)^s\le\exp(-s),
$$
whenever $k>N_4$ for some $N_4>0$ depending only on $\alpha$. 
Since $\alpha$ depends only on $\delta$, $N_3$ and $N_4$ also depend
only on $\delta$. Let $N=\max\{N_1,N_2,N_3,N_4\}$. Then $N$ depends
only on $\delta$ and provided $k>N$, the expected number of bad $S$
is at most
$$
 \sum_{\log^2 n\le s\le n/h}\exp(-s)=o(1).
$$
 Lemma~\ref{l2:SmallSet} follows. \qed
\ss

It only remains to prove Claim~\ref{c:p(q)}.\ss

\no {\bf Proof of Claim~\ref{c:p(q)}.\ } To illustrate the method of
computing $p(q,t)$, we show in detail the case $h=2$ first.
Conditional on that $d(S)=q$, we want to estimate the probability
that there are at least $t$ edges in $S$.
 Consider the alternative algorithm that generates the probability
 space of the
partition-allocation model $\mathcal{P}([n], [m_2],0,k+1)$. Fix any
allocation which allocates $q$ balls into bins representing vertices
in $S$ with each bin containing at least $k+1$ balls. There are at
most
$$
\binom{q}{2t}\frac{(2t)!}{2^tt!}
$$
partial partitions that contain $t$ parts within $S$.  The
probability of every such partial partition to occur is
$$
\prod_{i=0}^{t-1}\frac{1}{D-1-2i}.
$$
So
$$
p(q,t)\le
\binom{q}{2t}\frac{(2t)!}{2^tt!}\cdot\prod_{i=0}^{t-1}\frac{1}{D-1-2i},
$$
which is at most
\begin{eqnarray*}
&&\frac{[q]_t}{2^tt!}\prod_{i=0}^{t-1}\frac{q-t-i}{D-1-2i}\le
\left(\frac{eq}{2t}\right)^t\left(\frac{q-t}{D-1}\right)^t\le
\left(\frac{eq}{2t}\cdot\frac{q}{D}\right)^t
\end{eqnarray*}
Note that the second inequality holds since $q<D/2$ and so
$q-t<(D-1)/2$.

Now we estimate $p(q,t)$ in the general case $h\ge 2$. Let ${\bf
M}=([2m_{2}],\ldots,[hm_h])$. Consider the alternative algorithm
that generates the probability
 space of the
partition-allocation model $\mathcal{P}([n],{\bf M},{\bf 0},k+1)$,
defined in Section~\ref{core}. Fix any allocation that allocates
exactly $q$ balls into $S$ with each bin containing at least $k+1$
balls. The algorithm uniformly randomly partitions balls into parts
such that there are exactly $m_{h-j}$ parts with size $h-j$ for
$j=0,\ldots,w-1$. Let $\mathcal{U}=\{(u_2,\ldots,u_h)\in
I\hspace{-0.1cm} N^{(h-1)}:\sum_{i=2}^{h}(i-1)u_i=t\}$. Let ${\bf
u}=(u_2,\ldots,u_h)$ be an arbitrary vector from $\mathcal{U}$.
 We over estimate the probability that
$\rho(S,i)$ is at least $u_i$ for all $i=2,\ldots,h$, conditional on
$d(S)=q$. Let $p(q,{\bf u})$ denote this probability. Then clearly
$p(q,t)\le \sum_{{\bf u}\in \mathcal{U}}p(q,{\bf u})$. The number of
partial partitions that contain $u_i$ partial parts of size $i$
within $S$ is
\begin{equation}
\binom{q}{u_1,2u_2,3u_3,\ldots,hu_h}\frac{(2u_2)!}{2!^{u_2}u_2!}\cdots
\frac{(hu_h)!}{h!^{u_h}u_h!},\lab{eq:embeddings}
\end{equation}
where $u_1=q-\sum_{i=2}^h iu_i$. For any such partial partition we
compute the probability that it occurs. The algorithm starts from
picking a ball $v$ unpartitioned in $S$ and then it chooses at most
$h-1$
 balls that are
u.a.r.\ chosen from all the unpartitioned balls to be partitioned
into the part containing $v$.


The probability of the occurrence of a given $u_2$ partial parts of
size $2$ within $S$ is at most
$$
\prod_{i=0}^{u_2}(h-1)\frac{1}{D-1-hi}=(h-1)^{u_2}\frac{1}{D-1}\cdot\frac{1}{D-h-1}\cdots\frac{1}{D-1-h(u_2-1)}.
$$
The probability of the occurrence of a given $u_3$ partial parts of
size $3$ within $S$ is at most
\begin{eqnarray*}
\prod_{i=0}^{u_3-1}\binom{h-1}{2}\frac{1}{D-hu_2-hi-1}\cdot\frac{1}{D-hu_2-hi-2}\le
(h-1)^{2u_3}\prod_{i=0}^{u_3-1}\frac{1}{(D-hu_2-hi-1)^2}.
\end{eqnarray*}
Note that the above inequality holds because $h\sum_{i=2}^hu_i\le
hq/2<D/2$.
 Keeping the analysis in this procedure, we obtain that
the probability of a particular partial partition with $u_i$ partial
parts of size $i$ within $S$ is at most
\begin{eqnarray*}
&&(h-1)^{u_2+2u_3+\cdots+(h-1)u_h}\times\prod_{i=0}^{u_2-1}\frac{1}{D-hi-1}\prod_{i=0}^{u_3-1}\frac{1}{(D-hu_2-hi-1)^2}\\
&&{\hspace{1cm}}\times\cdots\times\prod_{i=0}^{u_{h-1}-1}\frac{1}{(D-h\sum_{j=2}^{h-2}u_j-hi-1)^{h-1}}
\prod_{i=0}^{u_{h-1}-1}\frac{1}{(D-h\sum_{j=2}^{h-1}u_j-hi-1)^{h-1}}.
\end{eqnarray*}
The product of this and~\eqn{eq:embeddings} gives an upper bound of
$p(q,{\bf u})$, which is at most
\begin{eqnarray}
&&\frac{[q]_{\sum_{i=2}^hiu_i}(h-1)^t}{u_2!u_3!\cdots
u_h!2!^{u_2}\cdots h!^{u_h}!}
\prod_{i=0}^{u_2-1}\frac{1}{D-hi-1}\cdots\prod_{i=0}^{u_h-1}\frac{1}{(D-h\sum_{j=2}^{h-1}u_j-hi-1)^{h-1}}.\nonumber
\end{eqnarray}
Since $2!^{u_2}\cdots h!^{u_h}!\ge 2^t$, this is at most
\begin{eqnarray}
&&\frac{[q]_t (h-1)^t}{u_2!u_3!\cdots
u_h!2^t}\prod_{i=0}^{u_2-1}\frac{q-t-i}{D-hi-1}\times\cdots\times\prod_{i=0}^{u_h-1}
\frac{q-t-\sum_{j=2}^{h-1}u_j-i}{(D-h\sum_{j=2}^{h-1}u_j-hi-1)^{h-1}}.\nonumber
\end{eqnarray}
Since $q<D/h$ and so $q-t\le (D-1)/h$, this is at most
\begin{eqnarray}
&&\frac{(eq(h-1)/2)^t}{u_2^{u_2}\cdots u_h^{u_h}}\left(\frac{q-t}{D-1}\right)^{u_2}
\left(\frac{q-t-u_2}{(D-hu_2-1)^2}\right)^{u_3}\cdots\left(\frac{q-t-\sum_{j=2}^{h-1}u_j}{(D-h\sum_{j=2}^{h-1}u_j-1)^{h-1}}\right)^{u_h}\nonumber\\
&&{\hspace{.5cm}}=\frac{(eq(h-1)/2)^t}{u_2^{u_2}\cdots u_h^{u_h}(q-t-u_2)^{u_3}\cdots(q-t-\sum_{j=2}^{h-1}u_j)^{(h-2)u_h}}\nonumber\\
&&{\hspace{1cm}}\times\left(\frac{q-t}{D-1}\right)^{u_2}\left(\frac{q-t-u_2}{D-hu_2-1}\right)^{2u_3}
\cdots\left(\frac{q-t-\sum_{j=2}^{h-1}u_j}{D-h\sum_{j=2}^{h-1}u_j-1}\right)^{(h-1)u_h}\nonumber\\
&&{\hspace{.5cm}}\le
\frac{(eq(h-1)/2)^t}{u_2^{u_2}(u_3(q-t-u_2))^{u_3}\cdots(u_h(q-t-\sum_{j=2}^{h-1}u_j)^{h-2})^{u_h}}\left(\frac{q-t}{D-1}\right)^{t}\nonumber\\
&&{\hspace{.5cm}}\le
\frac{(eq(h-1)/2)^t}{u_2^{u_2}(u_3(q-t-u_2))^{u_3}\cdots(u_h(q-t-\sum_{j=2}^{h-1}u_j)^{h-2})^{u_h}}\left(\frac{q}{D}\right)^{t}.\nonumber
\end{eqnarray}
 Since
$q\ge\sum_{i=2}^{h}iu_i$ and $t= \sum_{i=2}^{h}(i-1)u_i$,
$q-t-\sum_{j=2}^i u_j\ge \sum_{j=i+1}^h u_j\ge u_{i+1}$ for all
$2\le i\le h-1$, and so
$$
u_2^{u_2}(u_3(q-t-u_2))^{u_3}\cdots\Big(u_h\Big(q-t-\sum_{j=2}^{h-1}u_j\Big)^{h-2}\Big)^{u_h}\ge
u_2^{u_2}u_3^{2u_3}\cdots u_h^{(h-1)u_h}.
$$
We prove the following claim  later.
\begin{claim}\lab{c:opt} Let $t=\sum_{j=2}^h(j-1)u_j$. Then
$$
u_2^{u_2}u_3^{2u_3}\cdots u_h^{(h-1)u_h}\ge
\left(\frac{2t}{h(h-1)}\right)^t.
$$
\end{claim}
By Claim~\ref{c:opt}, for any $h\ge 2$,
$$
p(q,{\bf u})\le
\left(\frac{eqh(h-1)^2}{4t}\cdot\frac{q}{D}\right)^t,\ \ \ \forall
{\bf u}\in\mathcal{U}.
$$
Since $|\mathcal{U}|<t^h$, we have
$$
p(q,t)\le
t^h\left(\frac{eh(h-1)^2q^2}{4tD}\right)^t=\left(\exp\left(\frac{h\ln
t}{t}\right)\frac{eh(h-1)^2q^2}{4tD}\right)^t.
$$
In particular, if $t\rightarrow \infty$ as $n\rightarrow\infty$,
then $h\ln t/t\rightarrow 0 $ and so $\exp(h\ln t/t)\le (h/(h-1))^2$
provided $n$ is large enough. So
$$
p(q,t)\le\left(\frac{eh^3q^2}{4tD}\right)^t.\qed
$$

\no {\bf Proof of Claim~\ref{c:opt}.\ } We solve the following
optimization problem
\begin{eqnarray*}
(P_1)\ \  \min&&m_2^{m_2}m_3^{2m_3}\cdots m_h^{(h-1)m_h}\\
s.t.&& m_2+2m_3+\cdots+(h-1)m_h=t\\
&&m_2,m_3,\ldots,m_h\ge 0
\end{eqnarray*}

Letting $x_i=(i-1)m_i$ for $2\le i\le h$, and taking the logarithm
of the objective function, $(P_1)$ is equivalent to the following
optimization problem.
\begin{eqnarray*}
(P_2)\ \  \min&&x_2\ln x_2+x_3\ln (x_3/2)+\cdots +x_h \ln(x_h/(h-1))\\
s.t.&& x_2+x_3+\cdots+x_h=t\\
&&x_2,x_3,\ldots,x_h\ge 0
\end{eqnarray*}
For convention,  let $x\ln x=0$ if $x=0$. Applying the Lagrange
multiplier yields ${\bf x}^*=(x_2^*,x_3^*,\ldots,x_h^*)$ with
$x_i=2t(i-1)/h(h-1)$, which is a feasible solution of $(P_2)$. In
order to show that this is an optimal solution, we need to show that
the optimal solution does not appear on the boundary.

Let ${\bf x}$ be any solution on the boundary of $(P_2)$. Then there
exists $2\le i\le h$ such that $x_i=0$. There also exists $j$ with
$x_j>0$. Consider ${\bf x}'$ with $x'_i=(i-1)x_j/h$,
$x'_j=x_j-(i-1)x_j/h$ and $x'_l=x_l$ for any $l\neq i,j$. Then ${\bf
x}'$ is feasible and it is straightforward to check that
$$
x'_i\ln(x'_i/(i-1))+x'_j\ln(x'_j/(j-1))<x_i\ln(x_i/(i-1))+x_j\ln(x_j/(j-1)).
$$
Hence ${\bf x}'$ cannot be an optimal solution. This proves that
${\bf x}^*$ is the minimizer and so the optimal value of $(P_1)$ is
$\exp(t\ln(2t/(h(h-1))))=(2t/(h(h-1)))^t$.\qed \ss \remove{
 The objective function of $(P_2)$ is
$$
f(x_2,\ldots,x_h)=x_2\ln x_2+x_3\ln (x_3/2)+\cdots +x_h
\ln(x_h/(h-1)).
$$
So
$$
\frac{\partial}{\partial x_i}
f(x_2,\ldots,x_h)=1+\ln\left(\frac{x_i}{i-1}\right), \ \ \mbox{for}\
i=2,\ldots,h.
$$
Letting
$$
1+\ln\left(\frac{x_i}{i-1}\right)=y, \ \ \mbox{for all}\
i=2,\ldots,h,
$$
we have
$$
x_i=\frac{2t(i-1)}{h(h-1)}, \ \ \mbox{for all}\ i=2,\ldots,h.
$$
We prove that
$(x_2^*,x_3^*,\ldots,x_h^*)=\left(\frac{2t(i-1)}{h(h-1)}\right)_{i=2,\ldots
,h}$ is the optimal solution to $(P_2)$.

Any feasible solution of $(P_2)$ can be expressed as
$(x_2,x_3,\ldots,x_h)=(x_2^*+\Delta x_2,x_3^*+\Delta
x_3,\ldots,x_h^*+\Delta x_h)$, such that $x_2,x_3,\ldots,x_h\ge 0$,
and $\Delta x_2+\cdots +\Delta x_h=0$.
 Assume $(x_2,x_3,\ldots,x_h)$ is an optimal solution of $(P_2)$ and $(x_2,x_3,\ldots,x_h)\neq (x_2^*,x_3^*,\ldots,x_h^*)$.
 Then there must exists $i,j$, such that $\Delta x_i>0$ and $\Delta x_j<0$. Hence
$$
\frac{\partial}{\partial x_i}
f(x_2,\ldots,x_h)>\frac{\partial}{\partial x_j} f(x_2,\ldots,x_h)>0.
$$
Then there exists a small constant $\delta>0$, such that
$(x'_2,x'_3,\ldots,x'_h)=(x_2,\ldots,x_i-\delta,\ldots,x_j+\delta,\ldots,x_h)$
is a feasible solution to $(P_2)$, and
$f(x'_2,x'_3,\ldots,x'_h)<f(x_2,x_3,\ldots,x_h)$, contradicting
$(x_2,x_3,\ldots,x_h)$ being an optimal solution. Hence
$(x_2^*,x_3^*,\ldots,x_h^*)$ is the optimal solution. So the optimal
value of $(P_2)$ is
$$
f(x_2^*,x_3^*,\ldots,x_h^*)=t\ln\left(\frac{2t}{h(h-1)}\right).
$$
So the optimal value of $(P_1)$ is
$\exp(t\ln(2t/(h(h-1))))=(2t/(h(h-1)))^t$. Hence
$$
m_2^{m_2}m_3^{2m_3}\cdots m_h^{(h-1)m_h}\ge
\left(\frac{2t}{h(h-1)}\right)^t.\qed
$$

}


\end{document}